\documentclass[letterpaper,times]{sae}

\usepackage[svgnames]{xcolor}
\definecolor{SAEblue}{RGB}{1,160,233}
\raggedright

\usepackage[justification=justified, singlelinecheck=false]{caption}  
\DeclareCaptionFont{eightpt}{\fontsize{8pt}{11pt}\selectfont #1}
\usepackage{lineno}
\usepackage[utf8]{inputenc}

\usepackage{cite}

\usepackage{array}
\newcolumntype{L}[1]{>{\raggedright\let\newline\\\arraybackslash\hspace{0pt}}p{#1}}
\newcolumntype{C}[1]{>{\centering\let\newline\\\arraybackslash\hspace{0pt}}p{#1}}
\newcolumntype{R}[1]{>{\raggedleft\let\newline\\\arraybackslash\hspace{0pt}}p{#1}}

\usepackage{graphicx}
\usepackage{multirow}
\usepackage{amsmath}
\usepackage{amssymb}
\usepackage{comment}

\newcommand{\ignore}[1]{}

\usepackage{nameref}
\usepackage{booktabs}

\makeatletter
\def\@seccntformat#1{%
  \expandafter\csname c@#1\endcsname\c@section
  }
\makeatother

\usepackage{subfigure}
\usepackage{multimedia}
\usepackage{lipsum}

\PaperTitle{Concurrent Optimization of Vehicle Dynamics and Powertrain Operation Using Connectivity and Automation}
\usepackage{calc} 

\makeatletter 
\renewcommand\@biblabel[1]{#1. } 
\makeatother

\AddAuthor{A M Ishtiaque Mahbub, Andreas A. Malikopoulos}{University of Delaware, USA} 
\PaperNumber{20XX-XX-XXXX} 
\SAECopyright{2020} 

\begin{document}
\maketitle
\section{Abstract}
Connected and automated vehicles (CAVs) provide the most intriguing opportunity to reduce energy consumption and travel delays. In this paper, we propose a two-level control architecture for CAVs to optimize (1) the vehicle's speed profile, aimed at minimizing stop-and-go driving, and (2) the powertrain efficiency of the vehicle for the optimal speed profile derived in (1). The proposed hierarchical control framework can be implemented onboard the vehicle in real time with minimal computational effort. We evaluate the effectiveness of the efficiency of the proposed architecture through simulation in Mcity using a 100\% penetration rate of CAVs. The results show that the proposed approach yields significant benefits in terms of energy efficiency.

\section{Introduction} \label{sec:1}
\begin{figure}[h]
	\centering
	\includegraphics[width=3.5 in]{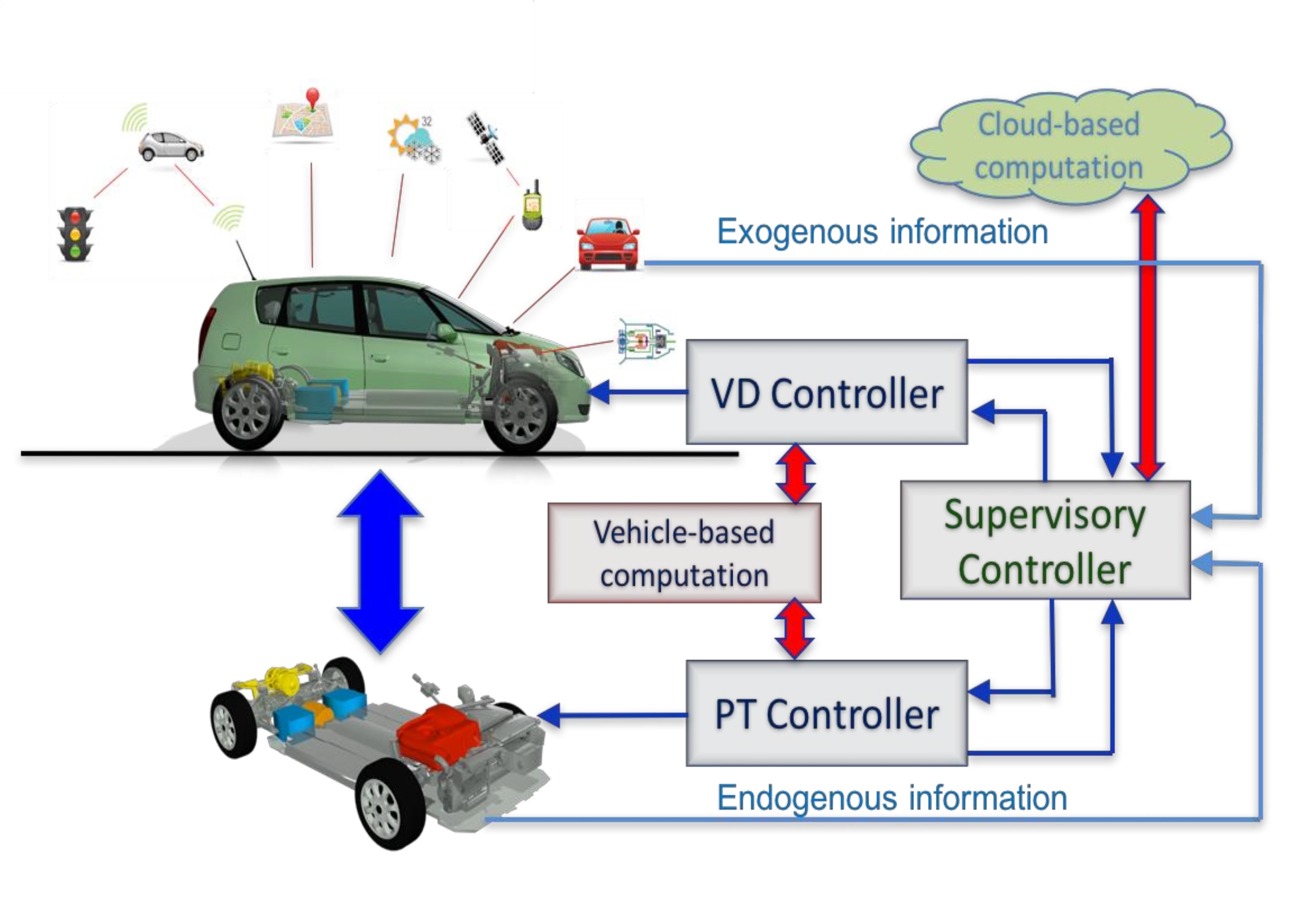} \caption{Vehicle's control architecture.}%
	\label{fig:supervisory-controller}%
\end{figure}

Recognition of the necessity for connecting vehicles to their surroundings has gained momentum. The main focus has been on safety and how accidents could be potentially prevented by developing multi-scale systems based on vehicle-to-vehicle (V2V) and vehicle-to-infrastructure (V2I) communications to alert drivers for a potential collision. The question is whether we could use connectivity and automation to optimize the efficiency of a vehicle in addition to safety. We are, in particular, interested in investigating the opportunities to improve the efficiency of hybrid electric vehicles (HEVs) and plug-in HEVs (PHEVs) \cite{Malikopoulos2014b} when these vehicles are connected and automated, and they can exchange information with their surrounding environment. Thus, we propose a two-level control architecture (Fig. \ref{fig:supervisory-controller}) for a PHEV. The objective is to (1) optimize the vehicle's speed profile aimed at minimizing stop-and-go driving, and (2) optimize the powertrain of the vehicle for this optimal speed profile obtained in (1). The control architecture consists of a vehicle dynamics (VD) controller, a powertrain (PT) controller, and a supervisory controller. The supervisory controller oversees the VD and PT controllers and communicate the endogenous and exogenous information appropriately. The VD controller optimizes online the acceleration/deceleration and speed profile of the vehicle in situations where there is a potential conflict with other vehicles, e.g., in traffic lights, stop signs, roundabouts, etc, so that to avoid stop-and-go driving. The PT controller computes the optimal nominal operation (set-points) for the engine and motor corresponding to the optimal solution of the VD controller. The complexity of the problem dimensionality can be managed by establishing two parallel and appropriately interacting computational levels: (1) a cloud-based, and (2) a vehicle-based level.

There is a solid body of research now available aimed at enhancing our understanding of optimizing the efficiency of both conventional vehicles \cite{Malikopoulos2008}  and HEVs or PHEVs \cite{Malikopoulos2014b}. Many different approaches have been proposed to address the fundamental vehicle system performance using offline and online analytical algorithms. It appears that research needs to be devoted to considering the vehicle as part of a larger system, which can be optimized at an even larger scale. Such large-scale optimization will require the acquisition and processing of additional information from the driver and surrounding environment. This is likely to require addition of new sensors and/or better utilization of information generated by existing sensors. However, the processing of such multi-scale information might require significantly new approaches in order to overcome the curse of dimensionality. Next, we attempt to summarize the work reported in the literature to date on optimizing vehicle-level and powertrain-level (specifically for HEVs and PHEVs) operations. Any such effort has obvious limitations. Space constraints limit the inclusion and detailed description of the rich literature in this area, and thus, we limited our efforts to reference work important for understanding the fundamental concepts or explaining significant departures from previous work. 

\subsection{Related Work}
\subsubsection{Vehicle Dynamics Optimization}
Intersections, merging roadways, roundabouts, and speed reduction zones along with the driver responses to various disturbances  \cite{Malikopoulos2013} are the primary sources of bottlenecks that contribute to traffic congestion. Connected and automated vehicles (CAVs) can be coordinated at these sources of bottlenecks to improve traffic flow. Several research efforts have been reported in the literature proposing either centralized or decentralized approaches on coordinating CAVs. One of the very early efforts in this direction was proposed in 1969 by Athans \cite{Athans1969} for safe and efficient coordination of merging maneuvers with the intention to avoid congestion. Assuming a given merging sequence, Athans  formulated the merging problem as a linear optimal regulator to control a single string of vehicles, with the aim of minimizing the speed errors that will affect the desired headway between each consecutive pair of vehicles. Varaiya \cite{Varaiya1993} discussed extensively the key features of an automated intelligent highway system. 

In 2004, Dresner and Stone \cite{Dresner2004} proposed the use of the reservation scheme to control a single-free intersection of two roads with vehicles traveling with similar speeds on a single direction on each road. Since then, several research efforts have extended this approach \cite{Dresner2008,DeLaFortelle2010, Huang2012}. Some approaches have focused on coordinating vehicles at intersections to improve the traffic flow \cite{Zohdy2012,Yan2009,kim2014}. Recently, a decentralized optimal control framework was established for coordinating online CAVs in different transportation segments. A closed-form, analytical solution was presented in \cite{Rios-Torres2} and \cite{Ntousakis:2016aa} for coordinating online CAVs at highway on-ramps, in \cite{Malikopoulos2017} and in \cite{Mahbub2019ACC} at intersections, in \cite{Malikopoulos2018a} at roundabouts, in \cite{Malikopoulos2018c} for speed harmonization, and in \cite{Zhao2018} in a corridor with several conflict zones. The solution of the optimal control problem considering state and control constraints was presented in \cite{Malikopoulos2017} and \cite{Mahbub2019CDC} at an urban intersection without considering rear-end collision avoidance constraint, and the conditions under which the latter does not become active were presented in \cite{Malikopoulos2018c}. These efforts considered a control zone inside of which the CAVs can communicate with each other through vehicle-to-vehicle (V2V) connectivity, and with a coordinator through vehicle-to-infrastructure (V2I) connectivity, giving rise to the decentralized communication topology. The coordinator is not involved in any decision making process. It only handles the communication of appropriate information among CAVs. The performance and the effectiveness of the aforementioned decentralized vehicle dynamics (VD) controllers for individual tasks (i.e., on-ramp merging, roundabout, signal-free intersection) have been validated at University of Delaware's 1:25 Scaled Smart City  \cite{beaver2019demonstration} using robotic cars. A recent survey paper \cite{Malikopoulos2016a} includes detailed discussions of the research reported in the literature to date on coordination of CAVs to improve vehicle-level operation.

\subsubsection{Vehicle Powertrain Optimization}
The power management control algorithm in HEVs and PHEVs determines how to split the power demanded by the driver between the thermal and electrical subsystems so that maximum fuel economy and minimum pollutant emissions can be achieved. Developing the control algorithm in HEVs and PHEVs constitutes a challenging control problem and has been the object of intense study for the last 20 years \cite{Malikopoulos2016AMO, Malikopoulos2014c,Malikopoulos2013a, Malikopoulos2015_ITS_HEV}. 
A significant amount of work has been proposed on optimizing the power management control in HEVs using dynamic programming (DP) \cite{Lin2003}. DP has been used to benchmark the fuel economy of HEVs by providing the maximum theoretical efficiency over a given driving cycle. It has been also extended to the stochastic formulation for a family of driving cycles \cite{Lin2004,Tate2010}. Some efforts have included the shortest path formulation \cite{Tate2007} for parallel HEV trucks. Opila \textit{et al.} \cite{Opila2008} presented a method to account for drivability metrics in their proposed power management control algorithm. Although DP can provide the optimal solution in both the deterministic and stochastic formulation of the power management control problem, the computational burden associated with deriving the optimal control policy prohibits online derivation in vehicles. To address these issues, research efforts have been concentrated on developing online algorithms consisting of an instantaneous optimization problem that accounts for storage system variation through the equivalent fuel consumption.  Paganelli \textit{et al.} \cite{Paganelli2001} introduced the equivalent consumption minimization strategy (ECMS) that optimizes the power split and the gear ratio while assigning a nonlinear penalty function for SOC deviation in a parallel HEV. Sciarretta, Back, and Guzzella \cite{Sciarretta2004} proposed an ECMS algorithm in which EFC is evaluated under the assumption that every variation in SOC will be compensated in the future by the engine running at the current operating point. Simulation results illustrated that the proposed algorithm can keep deviations of SOC from the target value at a low level. Musardo, Rizzoni, and Staccia \cite{Musardo2005} presented an adaptive ECMS (A-ECMS) algorithm that periodically computes the equivalence factor and refreshes the control parameters based on the current driving conditions to maximize fuel economy in a parallel HEV. In 2007, Pisu and Rizzoni \cite{Pisu2007} compared three algorithms that can be implemented online, a rule-based algorithm, an A-ECMS, and an $\mathcal{H}_\infty$ control. The simulation results showed that A-ECMS promises superior robustness and drivability, while it achieves better fuel economy results compared to the rule-based and $\mathcal{H}_\infty$ control algorithms. There has been also a significant amount of work using model predictive control but mainly in power split HEVs \cite{Borhan2012} and series HEVs (see \cite{Malikopoulos2013a} and the references therein). Some papers have focused on incorporating the destination routes \cite{Johannesson2007, Ambuhl2009}. To address variation in fuel consumption for different driving styles \cite{Malikopoulos2012,Malikopoulos2013d}, Huang, Tan, and He \cite{Huang2011} developed a statistical approach to distinguish automatically driving styles in HEVs. There have been also efforts to address the two vehicle-level and powertrain level operation simultaneously \cite{Luo2015a, Luo2015b, Li2017a}. However, these efforts have exhibited some limitation for online implementation. A detailed discussion of the research reported in the literature today on the power management control for HEVs/PHEVs can be found in \cite{Malikopoulos2014b}.
\subsubsection{Contribution and Structure}
In this paper, we propose a two-level control architecture (Fig. \ref{fig:supervisory-controller}) for a PHEV to (1) optimize the vehicle's speed profile aimed at minimizing stop-and-go driving and (2) optimize its powertrain efficiency for this optimal speed profile. We evaluate the effectiveness of the efficiency of the proposed architecture through simulation in a network of CAVs in Mcity, a 32 acre vehicle testing facility located at the north campus of University of Michigan. The contribution of this paper is the analysis of online optimization of the vehicle- and powertrain-level operation of a PHEV and classification of the improvements.

The remainder of the paper is organized as follows. First, we introduce the control architecture that consists of the VD controller and the PT controller. Then, we provide the simulation framework for the proposed VD and PT controller. Afterwards, we evaluate the effectiveness of the efficiency of the proposed approach in a simulation environment. Finally, we draw conclusions and discuss next steps.

\section{Control Architecture} \label{sec:2}
\begin{figure}[ht]
\centering
\includegraphics[width=2.5in]{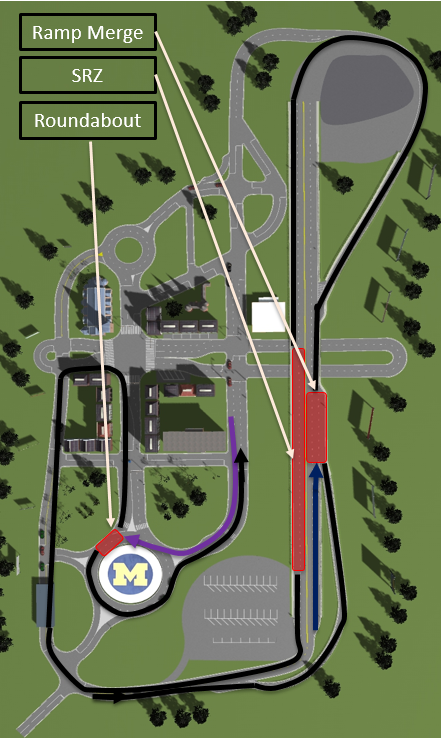} 
\caption{The corridor in Mcity with the conflict zones.}%
\label{fig:corridor}%
\end{figure}

\begin{figure}[h]
\centering
\includegraphics[width=3.5in]{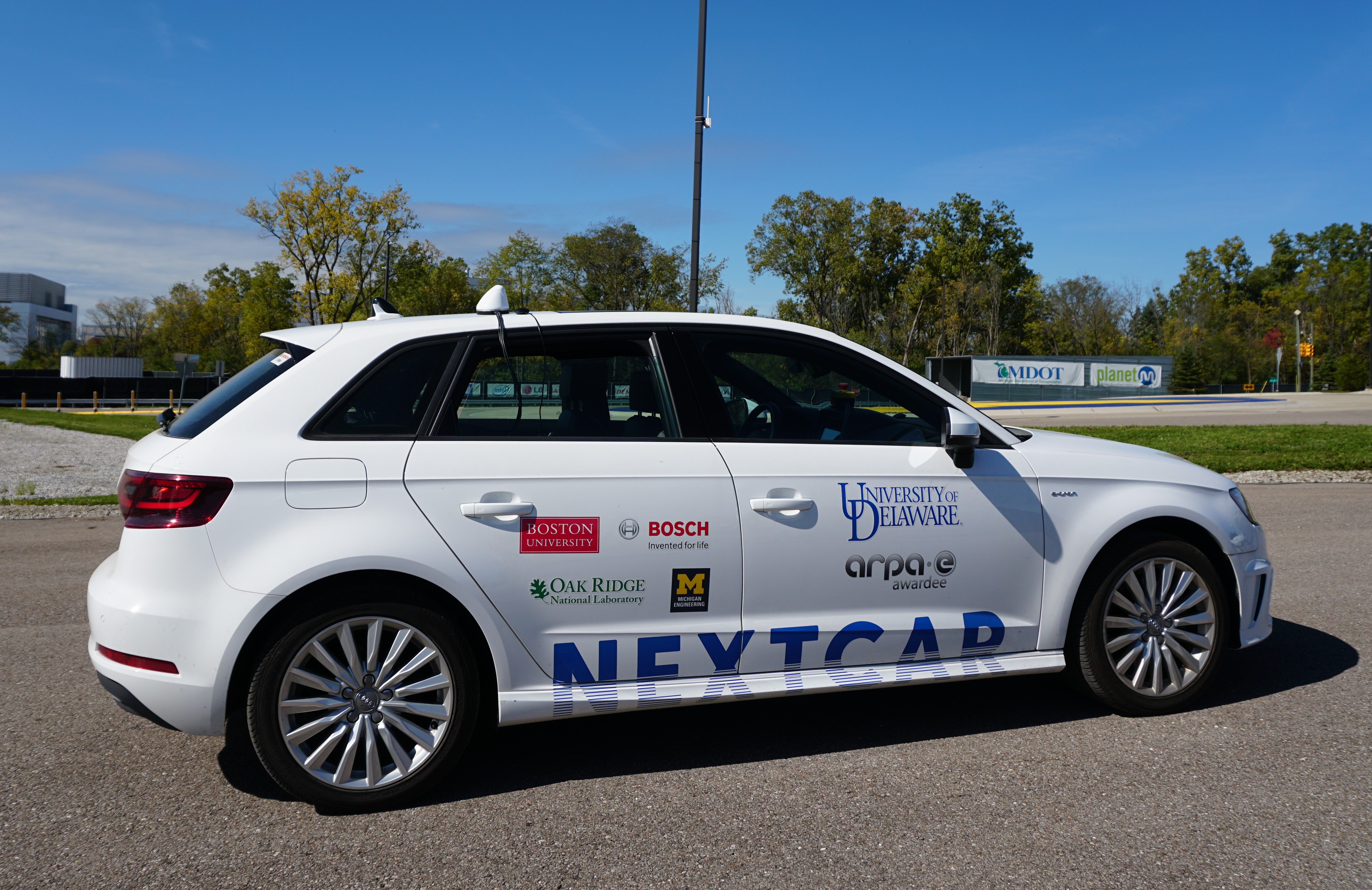} 
\caption{The target vehicle, Audi A3 e-tron, in Mcity.}%
\label{fig:audi}%
\end{figure}

We consider a network of connected and automated PHEVs (CA-PHEVs) driving through a corridor in Mcity  that consists of several conflict zones, e.g., a merging at roadways on-ramp, a speed reduction zone (SRZ), and roundabout as shown in Fig. \ref{fig:corridor}. The CA-PHEVs are retrofitted with necessary communication devices to interact with other CA-PHEVs and structures within their communication range through V2V and V2I communication. In this work, the CA-PHEV under consideration is an Audi A3 e-tron plug-in hybrid electric vehicle (Fig. \ref{fig:audi}).

The proposed control architecture applies to the operation of each CA-PHEV over a range of real-world driving scenarios deemed characteristic of typical commutes. A typical commute of a vehicle includes merging at roadways, crossing signalized intersections, cruising in congested traffic, and passing through speed reduction zones. Therefore, the VD controller needs to be able to optimize the speed profile of the vehicle over such traffic scenarios.  The PT controller computes the optimal nominal operation (set-points) for the engine and motor corresponding to the optimal speed profile provided by the VD controller. 

The focus on this paper is on the VD and PT controllers. The supervisory controller  coordinates the VD and PT controllers to ensure the optimal solution yielded by the VD controller is feasible for the PT controller and eventually results in maximization of the vehicle's energy efficiency. However, the details of this coordination along with the implications of the computational efforts of the VD and PT controllers are outside the scope of this paper.

We consider a corridor that contains three conflict zones $z = 1, 2, 3$ representing a merging roadway, a speed reduction zone, and a roundabout respectively. The corridor containing the CA-PHEV's main route is illustrated by the black trajectory in Fig. \ref{fig:corridor}. The highway, and the urban traffic routes are shown by the northbound blue and southbound purple trajectories respectively in Fig. \ref{fig:corridor}, and are termed as the \textit{congestion routes}. Vehicles travelling through the congestion routes create traffic congestion for the vehicles travelling through the main route. Upstream of each conflict zone, we define a \textit{control zone} (CZ) where the CA-PHEVs can coordinate with each other and pass the conflict zone without any rear-end or later collision. Each CZ has a coordinator that can communicate with the CA-PHEVs traveling within its range. Note that the coordinator does not get involved in any decision of the vehicles. The coordinator just assigns a unique identity to each CA-PHEV when they enter a CZ. We denote the time for each CA-PHEV $i$ to enter the CZ of the conflict zone $z$ to be $t_i^{0,z}$, the time to enter the conflict zone to be $t_i^z$, and the time to exit the CZ to be $t_i^{f,z}$. The objective of each CA-PHEV is to derive its optimal control input (acceleration) to cross each of the conflict zones without any rear-end or lateral collision with other vehicles, and simultaneously optimize its powertrain to achieve better energy efficiency. 

In the modeling framework described above, we impose the following assumptions:

\textbf{Assumption 1: }For each CA-PHEV $i$, none of the constraints is active when entering the CZ of conflict zone $z$ at time $t_i^{0,z}$.

\textbf{Assumption 2: }Each CA-PHEV is equipped with sensors to measure and share their local information while communication among CA-PHEVs occurs without any delays or errors. 

\textbf{Assumption 3: } No lane change is allowed within the vehicle's route.

Assumption 1 ensures that the initial state and control input are feasible. The second assumption may be strong, but it is relatively straightforward to relax as long as the noise in the measurements and/or delays is bounded. For example, we can determine upper bounds on the state uncertainties as a result of sensing or communication errors and delays, and incorporate these into more conservative safety constraints. The last assumption simplifies the modeling framework by restricting the vehicle route in a single lane.

Let $N(t)\in\mathbb{N}$ be the number of CA-PHEVs inside the CZ of a conflict zone $z= 1,2,3$ of at time $t\in\mathbb{R}^{+}$. We denote the sequence of the vehicles to be entering a CZ of conflict zone $z$ as $\mathcal{N}_z(t)=\{1,..., N(t)\}$. We use the double integrator model to represent the dynamics of each CA-PHEV $i$ as follows,
\begin{equation}%
\begin{split}
\dot{p}_{i} &  =v_{i}(t),\\
\dot{v}_{i} &  =u_{i}(t),
\label{eq:model2}
\end{split}
\end{equation}
where $p_{i}(t)\in\mathcal{P}_{i}$, $v_{i}(t)\in\mathcal{V}_{i}$, and $u_{i}(t)\in\mathcal{U}_{i}$ denote the position, speed and
acceleration/deceleration (control input) of each vehicle $i$ in the corridor. 
The sets $\mathcal{P}_{i}$, $\mathcal{V}_{i}$, and $\mathcal{U}_{i}$ are complete and totally bounded subsets of $\mathbb{R}$.
Let $x_{i}(t)=\left[p_{i}(t) ~ v_{i}(t)\right] ^{T}$ denote the state of each vehicle $i$, with initial value $x_i(t_i^{0,z})= x_{i}^{0,z}=\left[p_{i}^{0,z} ~ v_{i}^{0,z}\right]^{T}$ taking values in $\mathcal{X}_{i}%
=\mathcal{P}_{i}\times\mathcal{V}_{i}$. Here, $p_{i}^{0,z} = p_i(t_i^{0,z}) $, and $ v_{i}^{0,z} = v_i(t_i^{0,z})$. The state space $\mathcal{X}_{i}$ for each vehicle $i$ is closed with respect to the induced topology on $\mathcal{P}_{i}\times \mathcal{V}_{i}$ and thus, it is compact.

To ensure that the control input and vehicle speed are within a
given admissible range, the following constraints are imposed.
\begin{gather}%
u_{min} \leq u_{i}(t)\leq u_{max},\quad\text{and}\\
0  \leq v_{min}\leq v_{i}(t)\leq v_{max},\quad\forall t\in\lbrack t_{i}%
^{0,z},t_{i}^{f,z}],
\label{eq:speed_accel constraints}%
\end{gather}
where $u_{min}$, $u_{max}$ are the minimum deceleration and maximum
acceleration, and $v_{min}$, $v_{max}$ are the minimum and maximum speed limits respectively. Note that, we assume homogeneity in terms of vehicle types, which enables the use of constant maximum and minimum acceleration values for any vehicle $i$.

To ensure the absence of rear-end collision of two consecutive vehicles traveling on the same lane, the position of the preceding CA-PHEV $k$ should be greater than or equal to the position of the following vehicle $i$ plus a predefined safe distance $\delta_i(t)$. Thus we impose the rear-end safety constraint 
\begin{equation}
s_{i}(t)=p_{k}(t)-p_{i}(t) \geq \delta_i(t),~ \forall t\in [t_i^{0,z}, t_i^{f,z}],
\label{eq:rearend_constraint}
\end{equation}
where $s_{i}(t)$ denotes the distance of vehicle $i$ from vehicle $k$ which is physically immediately ahead of $i$. Relate the minimum safe distance $\delta_i(t)$ as a function of speed $v_i(t)$, 
\begin{equation}
\begin{split}
\delta_i(t)=\gamma_i + \rho_i \cdot v_i(t),~ \forall t\in [t_i^{0,z}, t_i^{f,z}],
\label{eq:safedist}
\end{split}
\end{equation}
where $\gamma_i$ is the standstill distance, and $\rho_i$ is minimum time gap that vehicle $i$ would maintain while following another vehicle.

\subsection{Vehicle Dynamics (VD) Controller}

When a vehicle enters the CZ, the coordinator receives its information and assigns a unique identity $i$ to the vehicle. 
The order of the vehicle $i \in \mathcal{N}_z(t)$ has to satisfy the following condition,
\begin{equation}\label{eq:MS}
t_{i}^{0,z} \ge t_{i-1}^{0,z}, ~ \forall i\in \mathcal{N}_z(t), \, i>1,
\end{equation}
where $t_{i-1}^{0,z}$ is the time that vehicle $(i-1)$ will be entering the CZ of conflict zone $z$.

The sequence $\mathcal{N}_z(t)$ will remain unchanged if Eq. \eqref{eq:MS} holds. Otherwise, the sequence $\mathcal{N}_z(t)$ will be updated by some other vehicle coordination scheme to change the order of the CA-PHEVs entering a CZ. In this work, we adopt the first-in-first-out queue to generate the sequence. To satisfy the rear-end safety constraint in Eq. \eqref{eq:rearend_constraint} at $t_i^z$, we impose the following condition.

\begin{equation}
    t_i^{z} = max\left \{min\left\{t_{i-1}^{z}+\frac{\delta(v_i(t))}{v_{i-1}(t_{i-1}^{z})}, \frac{L_z}{v_{min}}\right \}, \frac{L_z}{v_0(t_i^{0,z})}, \frac{L_z}{v_{max}}\right \},
\end{equation}

where $\delta(v_i(t))$ is the safety distance based on the location of the previous vehicle $(i-1)$ with respect to the vehicle $i$, and $L_z$ is the length of the CZ length of the conflict zone $z$.

For each vehicle $i\in\mathcal{N}_z(t)$ traveling towards a conflict zone $z= 1,2, 3$, we define the cost function $J_{i}(u(t))$ in $[t_i^{z,0}, t_i^{z}]$
\begin{gather}\label{eq:decentral_general}
J_{i}(u(t))=  \int_{t_i^{z,0}}^{t_i^{z}} C_i(u_i(t))~dt, ~\forall z=1,2, 3,\\ 
\text{subject to}%
:\eqref{eq:model2},\eqref{eq:speed_accel constraints},\eqref{eq:rearend_constraint},\text{
}p_{i}(t_i^{z,0})=0\text{, }p_{i}(t_i^{z})=p_z,\nonumber\\
\text{and given }t_i^{z,0}\text{, }v_{i}(t_i^{z,0})\text{, }t_i^{z}.\nonumber
\end{gather}
Here, $C_i(u_i(t))$ is a function of the control input $u_i(t)$ and can be viewed as a measure of energy. When $C_i(u_i(t))$ is considered as the $L^2$-norm of the control input, i.e. $C_i(u_i(t))=\frac{1}{2}u_i^2(t)$, the transient engine operation can be minimized which, eventually, represents the minimization of energy consumption \cite{Zhao2019CCTAa}. The solution of Eq. \eqref{eq:decentral_general} considering the state and control constraints is recursive, and quite involved in nature. To derive the analytical solution of Eq. \eqref{eq:decentral_general} with state and control constraints, we adopt the standard methodology used in optimal control problems with interior point constraints. If the initially derived unconstrained solution violates any of the state or control constraints, then the unconstrained arc is pieced together with the arc corresponding to the violated constraint. We then re-solve the problem with the two arcs pieced together. The two arcs yield a set of algebraic equations which are solved simultaneously using the interior and terminal boundary conditions. If the resulting solution containing the switching between the two arcs violates yet another constraint, then the last two arcs are pieced together with the arc responsible for the latest constraint violation, and we re-solve the problem with the three arcs pieced together. This process is repeated until the solution does not violate any other constraints. Due to space limitation, we omit the derivation of the analytic solution of Eq. \eqref{eq:decentral_general}.
The solution of the constrained problem has been addressed in \cite{Malikopoulos2017}, and it requires the constrained and unconstrained arcs of the state and control input to be pieced together as discussed above to satisfy the Euler-Lagrange equations and necessary condition of optimality.

\subsection{Powertrain (PT) Controller}
The CA-PHEV considered here has a parallel configuration, where the gasoline engine and the electric motor can to provide necessary power to the wheel either independently or in combination. The engine, which can be fully decoupled in electric motor only operation, is connected to the integrated motor generator (IMG) unit through a singular clutch, which is in turn connected to a dual clutch transmission. The electric motor is coupled with the engine and gearbox, and can act as a generator for charging the battery.


Following the modeling framework in \cite{Malikopoulos2016AMO}, the evolution of the CA-PHEV state is modeled as a controlled Markov chain with a finite state space,  $\mathcal{S} \subset \mathbb{R}^n$, and a finite control space, $\mathcal{U} \subset \mathbb{R}^m, n,m\in\mathbb{N}$, from which the power management controller selects control actions. In our formulation the state space is the entire range of the engine and motor speed, $\mathcal{S} \subset \mathbb{R}^2$, where the engine and motor speed progress in a compact subset of $\mathbb{R}$. The control space $\mathcal{U}$ is the vector of engine and motor torque, $\mathcal{U} \subset \mathbb{R}^2$.

The evolution of the state occurs at each of a sequence of stages $t=0,1,...$, and it is portrayed by the sequence of the random variables  $X_{t(1:2)}=(X_{t(1)},  X_{t(2)}))^T=(N_{eng},  N_{mot})^T\in\mathcal{S}$ and $U_{t(1:2)}=(U_{t(1)},  U_{t(2)}))^T=(T_{eng},  T_{mot})^T\in\mathcal{U}$, corresponding to the HEV state (engine and motor speed) and control action (engine torque and motor torque) respectively. For each state $X_{t(1:2)}=i\in\mathcal{S}$ a nonempty set $\mathcal{C}(i)\subset\mathcal{U}$ of admissible control actions (engine and motor torque) is given which implies that at each state $i\in\mathcal{S}$, the control action set $\mathcal{C}(i)\subset\mathcal{U}$ should include only the control actions that satisfy the physical constraints of the engine and the motor.

At each stage $t$, the controller observes the engine and motor speed,  $X_{t(1:2)}=i\in\mathcal{S}$, which is a function of the vehicle speed, and executes an action, $U_{t(1:2)}=\mu(X_{t(1:2)})$ (engine and motor torque), from the feasible set of actions, $U_{t(1:2)}\in\mathcal{C}(i)$, at that state. At the same stage $t$, an uncertainty, $W_{t(1:2)}$, is incorporated in the system consisting of the torque demanded by the driver as designated by the pedal position, e.g., accelerator or brake. At the next stage, $t+1$, the system transits to the state $X_{t+1(1:2)} =j\in\mathcal{S}$ and a one-stage expected cost, $k(X_{t(1:2)},U_{t(1:2)})$, is incurred corresponding to the engine's and motor's efficiency. After the transition to the next state, a new action is selected and the process is repeated. The state transition from one state to another is imposed by a discrete-time equation that describes the dynamics of the CA-PHEV. 

The objective of the PT controller is to derive a  control policy that minimizes the long-run expected average cost of the CA-PHEV to split torque demanded by the driver between the engine and the motor for the optimal speed profile derived by the VD controller as a solution of Eq. \eqref{eq:decentral_general}. For the power management control problem here, we select the average cost criterion as we wish to optimize the efficiency of each CA-PHEV on average, hence

\begin{equation}
J^\pi= \lim_{T\to\infty}\frac{1}{T+1}\mathbb{E}^\pi \left [\sum_{t=0}^{T} k(X_{t(1:2)},U_{t(1:2)})\right],\label{eq:2}
\end{equation}
where $k(X_{t(1:2)},U_{t(1:2)})$ is the one-stage cost of CA-PHEV.

However, the computational burden associated with deriving the optimal control policy in Eq. \eqref{eq:2} prohibits online derivation onboard a vehicle. It has been shown \cite{Malikopoulos2016b} that the optimal control policy in Eq. \eqref{eq:2} is equivalent to the Pareto control policy that can be derived by formulating a multiobjective problem. The latter consists of the engine's efficiency, $\eta_{eng}$, and the motor's efficiency, $\eta_{mot}$.  Given the engine and motor speed $X_{t(1:2)}$, the objective is to find the optimal control action $U_{t(1:2)}$ (engine and motor torque) that optimizes a multiobjective function reflecting both the engine's and the motor's efficiency. Hence, one of the objectives is the engine's efficiency

\begin{equation}
f_1(N_{eng}, T_{eng}) = \eta_{eng} ,
\label{eq:6}
\end{equation}
and the other one is the motor's efficiency, 

\begin{equation}
f_2(N_{mot},T_{mot}) = \eta_{mot}.
\label{eq:7}
\end{equation}
The multiobjective optimization problem is formulated as
\begin{gather} 
\min_{U_t} k(X_{t(1:2)}, U_{t(1:2)})= \nonumber\\
\max_{U_t} \big(\alpha\cdot f_1(X_{t(1)}, U_{t(1)}) + (1-\alpha)\cdot f_2(X_{t(2)}, U_{t(2)})\big),\label{eq:objective_function} \\
\text{s.t.} \ \sum_{i=1}^{2}U_{t(i)}= T_{driver},\nonumber
\end{gather}
where $\alpha$ is a scalar that takes values in [0,1], $X_{t(1:2)}=(X_{t(1)},  X_{t(2)}))^T=(N_{eng},  N_{mot})^T\in\mathcal{S}$, $U_{t(1:2)}=(U_{t(1)},  U_{t(2)}))^T=(T_{eng},  T_{mot})^T\in\mathcal{U}$ is the vector of engine and motor torque. The multiobjective optimization problem in Eq. \eqref{eq:objective_function} yields the Pareto efficiency set between the engine and the motor by varying $\alpha$ from $0$ to $1$ at any given state of the HEV. 
For each state of the CA-PHEV and torque demand, we derive the Pareto efficiency set of Eq. \eqref{eq:objective_function} offline, and store it in a table. If there are multiple solutions, then one of these solutions is selected randomly since all of them will yield the same one-stage expected cost. The Pareto control policy is then implemented online using this table. 


\section{Simulation Environment}
\subsection{VISSIM Network}
To evaluate the effectiveness of the proposed control architecture, we design a simulation scenario using the commercial software PTV VISSIM \cite{ptv2014ptv}. We create a simulation environment resembling the Mcity vehcile testing facility, and define the network routes for all the vehicles. In terms of the nature of vehicle control, we consider two scenarios: 
\begin{enumerate}
    \item \textbf{Baseline Vehicle Dynamics:} All vehicles are human-driven without any connectivity and automation. The VISSIM employs the \textit{Wiedemann} car following model \cite{wiedemann1974} to emulate the behavior of the human driven vehicles. However, the VISSIM built-in car following models are slightly different from the original Wiedemann model in \cite{wiedemann1974} due to the inclusion of certain parameters to introduce additional randomness and heterogeneity in terms of driving behavior. We apply VISSIM's Wiedemann-74 car following model for the urban traffic, and Wiedemann-99 for the freeway traffic. The conflict zones inside the corridor have conventional traffic signals, which the vehicles must abide by. We model the traffic signals (yield behavior) at the on-ramp merging, and the roundabout by imposing VISSIM's \textit{priority rule} object in the conflict areas. To model the speed reduction zone, we apply VISSIM's \textit{reduced speed area} object with specified route length and speed.
    
    \item \textbf{Optimal Vehicle Dynamics:} In this case, we have CA-PHEVs (i.e., connected and automated PHEVs) travelling through the corridor. The CA-PHEVs employ the VD controller to optimize its speed profile for increasing fuel efficiency. We consider automated conflict zone, where the conventional traffic signals are not present. We modify one of VISSIM's API, namely the DriverModel.dll to implement the VD controller written in C++ for the CA-PHEVs in the optimal controlled case. Within each CZ, we override each vehicle's internal car-following model with the DriverModel.dll containing the VD controller logic. Outside the CZ, the vehicles can switch back to the VISSIM's built-in Wiedemann car following model. At each simulation time step, each vehicle access the external DriverModel.dll, and computes the optimal control output based on it's location in the route.
\end{enumerate}

To investigate the robustness of the VD controller through different conflict zones of the corridor, we consider three different traffic volumes in the VISSIM's traffic network. The traffic flow in both the main route and the congestion routes are modified to achieve different traffic scenarios. Table \ref{tab:traffic_cases} presents the different traffic flows considered to achieve the low, medium, and high traffic volume.
\begin{table}[!htb]
\fontsize{8}{10}\selectfont
\centering
\caption{Traffic volume at different routes within Mcity network.}\label{tab:traffic_cases}
\begin{tabular}{| L{0.4\columnwidth-2\tabcolsep-1.2\arrayrulewidth} | L{0.2\columnwidth-2\tabcolsep-1.2\arrayrulewidth} | L{0.2\columnwidth-2\tabcolsep-1.2\arrayrulewidth} | L{0.2\columnwidth-2\tabcolsep-1.2\arrayrulewidth} | 
}
\hline
           \textbf{} & \textbf{Low}         & \textbf{Medium}         & \textbf{High}                      \\ \hline
\textbf{Main Route [vph/lane]} & 500      & 400   & 300           \\  \hline
\textbf{Highway [vph/lane]} & 800     &  600    & 400          \\  \hline
\textbf{SRZ [vph/lane]} & 1400     &  1100    & 800         \\  \hline
\textbf{Roundabout [vph/lane]} & 700    &  550   & 400             \\  \hline
	\end{tabular}
	\par
   \vspace{-0.15\skip\footins}
   \renewcommand{\footnoterule}{}
\end{table}

The main route has a length of 1500 $m$ inside the corridor in Mcity (Fig. \ref{fig:corridor}). 
The length of the control zone for the on-ramp,  the SRZ and the roundabout were selected to be 100 $m$. The parameters relevant to the VISSIM simulation environment are compiled in Table \ref{tab:sim_param}.

\begin{table}[!htb]
\fontsize{8}{10}\selectfont
\centering
\caption{Simulation Parameters}\label{tab:sim_param}
\begin{tabular}{| L{0.75\columnwidth-2\tabcolsep-1.2\arrayrulewidth} | L{0.25\columnwidth-2\tabcolsep-1.2\arrayrulewidth} |  }
\hline
\textbf{Vehicle Parameters} &                               \\ \hline
{Maximum Acceleration [$m/s^2$]} & 1.5                 \\  \hline
{Maximum Deceleration [$m/s^2$] }    & 3.0                   \\  \hline
{Safe Time Headway [$s$] }    & 1.2                   \\  \hline
\textbf{Traffic Network} &                               \\ \hline
{Corridor Length [$m$] }    & 1500                   \\  \hline
{Control Zone Length [$m$]} & 100                  \\  \hline
{SRZ Length [$m$]} & 125                  \\  \hline
\textbf{Speed Limit}    &                    \\  \hline
{On-Ramp Merging [$mph$] }    & 40                   \\  \hline
{Speed Reduction Zone [$mph$] }    & 18.6                   \\  \hline
{Roundabout [$mph$] }    & 25                   \\  \hline
	\end{tabular}
	\par
   \vspace{-0.15\skip\footins}
   \renewcommand{\footnoterule}{}
\end{table}

\subsection{Powertrain Model}
\begin{figure}[ht]
	\centering
	\includegraphics[width=3.5in]{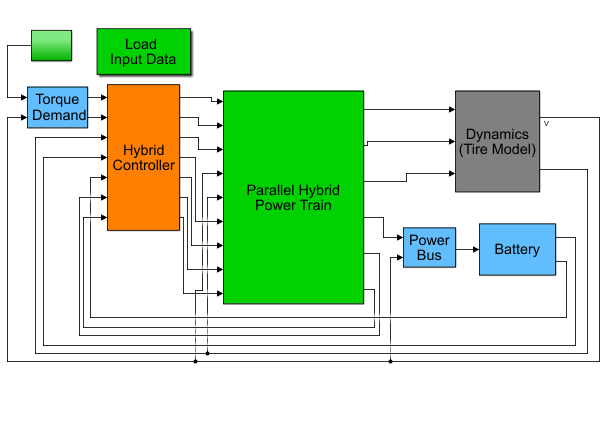} \caption{VESIM Model for modelling the powertrain of the plug-in hybrid electric vehicle.}%
	\label{fig:vesim}%
\end{figure}

To model the PT controller, we adopt a hybrid electric vehicle simulation tool VESIM \cite{vesim-1000037126/PAPER}, and create separate powertrain models for the baseline scenario and the optimized scenario. The purpose of this model is to represent the hybrid PT controller architecture of an HEV. The general architecture of the VESIM model is illustrated in Fig. \ref{fig:vesim}. We modify the VESIM to generate two different vehicle powertrain models as follows.

\subsubsection{Baseline VESIM (Audi A3 Powertrain)}
The baseline VESIM model has been calibrated to reproduce the vehicle characteristics generated by the factory controller of the Audi A3. Due to the combined contribution of the internal combustion engine and the motor of the Audi A3, the VESIM model computes the \textit{miles-per-gallon of gasoline equivalent} (MPGe) according to the EPA standard.
By feeding the baseline and the VD controller's speed profiles to the VESIM model, we quantify the fuel consumption of the Audi A3, and evaluate the performance of the VD controller at different conflict scenarios. Some of the most essential VESIM parameters required for the calibration purpose are summarized in Table \ref{tab:vesim_param}. 

The Audi A3 has several modes of operation:
\begin{enumerate}
    \item EV mode: Motor only operation, where the engine remains turned off.
    \item Charge Battery Mode: The engine provides the torque demanded by the driver, and also provides torque to the IMG unit to charge the battery.
    \item Hold Battery Mode: The SOC of the battery is maintained within a certain bandwidth of the initial SOC. Although the engine provides the torque demanded, but the motor can contribute if the torque demanded is higher than maximum capacity of the engine. 
    \item Hybrid Mode: The vehicle can use both the engine and the IMG unit to provide the torque demanded by the driver.
\end{enumerate}

\begin{figure}[ht]
	\centering
	\includegraphics[width=3.5in]{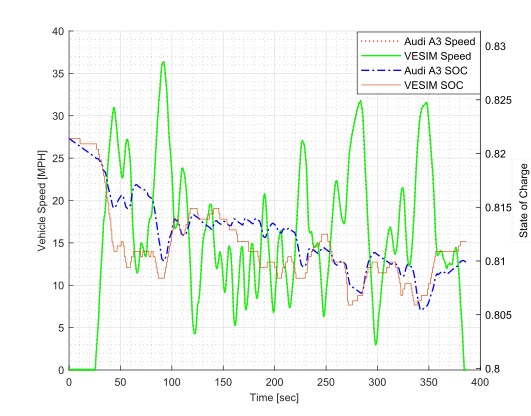} \caption{VESIM calibration model for hold battery mode.}%
	\label{fig:vesim_calibration}%
\end{figure}

We calibrate the VESIM model for each of the aforementioned modes of Audi A3 to capture appropriate engine and IMG characteristics. Fig. \ref{fig:vesim} shows a representation of such effort, where the speed and battery SOC profile of the Audi A3 and corresponding calibrated VESIM model for hold battery mode is illustrated. Note that, the hold battery mode is characterized by the SOC variation constrained within a certain SOC bandwidth.
We observe that the calibrated VESIM model can trace the reference speed of the actual Audi A3's drive cycle very closely, which is illustrated by their complete overlap in Fig. \ref{fig:vesim_calibration}. We obtain MPGe of 29.73 from the calibrated VESIM model compared to the Audi A3 MPGe of 30.92 for the same drive cycle, with 3.8 \% deviation.

\begin{table}[!htb]
\fontsize{8}{10}\selectfont
\centering
\caption{VESIM parameters calibrated with Audi A3 specification.}\label{tab:vesim_param}
\begin{tabular}{| L{0.65\columnwidth-2\tabcolsep-1.2\arrayrulewidth} | L{0.35\columnwidth-2\tabcolsep-1.2\arrayrulewidth} |  }
\hline
\textbf{Vehicle Parameters}    &                    \\  \hline

{Tire Model}    &      Michelline 225/60 r16              \\  \hline
{Weight (No driver) [$lb$]} & 3616 \\ \hline
{Rolling Resistance Coeff.} & 0.010 \\ \hline
{Frontal Area [$in$ x $in$]} & 56.1 x 60 \\ \hline
{Traction Torque Loss [\%]} & 0.95 \\ \hline
{Aerodynamic Drag Coeff.} & 0.32 \\ \hline
{Tire Traction Efficiency} & 0.96 \\ \hline
{Maximum Braking Force [$N$]} & 12000 \\ \hline

\textbf{Transmission}    &                    \\  \hline
{Gear Ratio [1-6]} & 3.50, 2.77, 1.85, 1.02, 1.02, 0.84 \\ \hline
{Gear Efficiency [1-6]} & 0.98, 0.98, 0.98, 0.98, 0.98, 0.98 \\ \hline
{Gear Intertia [1-6]} & 0.0023,	0.0009,	0.0023,	0.0009,	0.0023,	0.0009 \\ \hline
{Forward Drive Ratio} & 3.75 \\ \hline
{Forward Drive Efficiency} & 0.966 \\ \hline

\textbf{Engine (TFSI) Parameters} &                               \\ \hline
{Cylinder Volume [$cc$]}     &      1395              \\  \hline
{Maximum Engine Power [$kW$]} & 110                 \\  \hline
{Maximum Engine Torque [$Nm$]} & 250                 \\  \hline
{Engine Speed at Peak Torque [$rpm$]} & 1750~4000                 \\  \hline
{Engine Inertia [$kg/m^2$]} & 0.15                 \\  \hline

\textbf{Battery Pack} &                               \\ \hline
{Capacity [$kWh$] }    & 8.8                   \\  \hline
{Number of Cell/Module } & 12                 \\  \hline
{Number of Modules} & 8                  \\  \hline
{Maximum Voltage/Cell [$volt$]} & 4.2                  \\  \hline
{Minimum Voltage/Cell [$volt$]} & 2.1                  \\  \hline
{Maximum Battery Power [kW]} & 75                  \\  \hline

\textbf{IMG Unit}    &                    \\  \hline
{Maximum Motor Power [$kW$] }    & 100                   \\  \hline
{Maximum Motor Torque [$Nm$]} & 300                 \\  \hline
{Maximum Generator Torque [$Nm$]} & -300                 \\  \hline
{Motor Speed at Peak Torque [$rpm$]} & 2000                 \\  \hline
{Rotor's Rotational Inertia [$kg/m^2$]}   & 0.1 \\ \hline

\textbf{Energy Conversion (MPGe)}    &                    \\  \hline
{Gallon to Equivalent $CO_2$ }    & 8.887e-3                   \\  \hline
{kWh to Equivalent $CO_2$ }    & 7.44e-4                   \\  \hline
	\end{tabular}
	\par
   \vspace{-0.15\skip\footins}
   \renewcommand{\footnoterule}{}
\end{table}

\subsubsection{VESIM With PT Controller }

\begin{figure}[ht]
	\centering
	\includegraphics[width=3.5in]{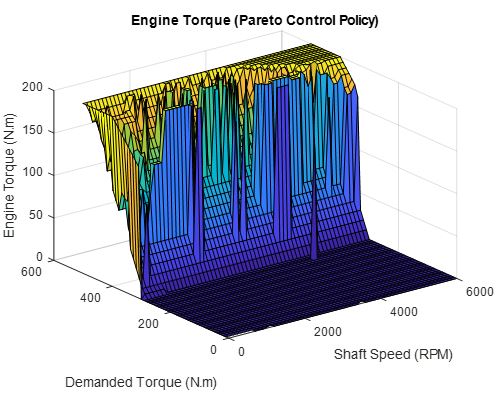} \caption{Pareto efficiency set of the powertrain controller.}%
	\label{fig:pareto_result}%
\end{figure}

For deriving the optimal control policy of the PT controller, we first derive the engine's efficiency map from its brake specific fuel consumption (BSFC) data, contrary to the motor's efficiency map which is readily available. With the engine's and motor's efficiency map, we solve the multiobjective optimization problem in Eq. \eqref{eq:objective_function} offline. We discretize the torque, engine/motor speed and the scalar $\alpha$ with the resolution of 10 $Nm$, 100 $RPM$ and 0.05 respectively, and use the optimization toolbox of MATLAB to solve Eq. \eqref{eq:objective_function}. Finer resolution of $\alpha$ increases the solution time significantly without yielding proportional changes to the solution matrix. The optimization process yields a Pareto efficiency set that we store in the CA-PHEV memory for later online use. For each torque demanded by the driver and corresponding engine/motor speed, the CA-PHEV searches the Pareto efficiency table to obtain the optimal torque split between the engine and the IMG unit.

The Pareto efficiency computed off-line is illustrated in Fig. \ref{fig:pareto_result}.
We note that when the driver's torque demand is below 300 $Nm$, the optimal solution is to use the motor exclusively to satisfy the torque demanded by the driver. This is because the electric motor considered here has high efficiency (almost 95\%) in most of its operating region compared to the engine, which has a peak efficiency at 35\%.

\section{Results and Discussion}

\subsection{PT Controller Performance}
The impact of the PT controller on engine operation is shown in Fig. \ref{fig:operating-points}. For the baseline scenario, we use the calibrated VESIM model, while the PT controlled CA-PHEVs were operated by the optimal VESIM model embedded with the Pareto efficiency table. We first evaluate the impact of the PT controller under different driving behavior. To ths end, we use three standardized drive cycles, namely the highway fuel economy driving schedule (HWFET), urban dynamometer driving schedule (UDDS), and the US06 supplemental federal test procedure to represent the 60mph-highway, heavy duty urban, and high acceleration aggressive driving behavior. We incorporate the aforementioned driving behaviors in a single drive cycle by stitching the considered drive cycles together to obtain a combined cycle of 25.72 miles. We characterize the performance of the proposed Pareto efficient powetrain controller compared to the baseline Audi A3 powertrain by tracing the drive cycles through the corresponding VESIM model.
\begin{figure}[htb]
	\centering
	\includegraphics[width=3.5in]{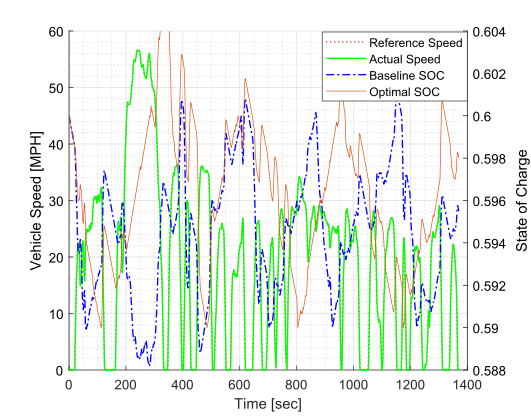}\caption{The UDDS drive cycle traced by the baseline (Audi A3) powertrain and the optimal controlled Pareto efficient powertrain. }%
	\label{fig:udds}%
\end{figure}
Fig. \ref{fig:udds} shows the UDDS drive cycle traced by the VESIM for the baseline (Audi A3 powertrain) and the optimal controlled (Pareto efficient) case. Note that, in both cases, the VESIM model was able to trace the UDDS drive cycle completely, as represented by the complete overlap of the reference and actual vehicle speed profile in Fig. \ref{fig:udds}. The battery SOC for the baseline and the optimal PT controlled case are also illustrated at the right axis of Fig. \ref{fig:udds}. We observe that, the battery SOC is constrained within a certain bandwidth of the initial SOC to represent the Audi A3's hold-battery mode.

The energy consumption results (MPGe) of the Audi A3's baseline powertrain and the optimal PT controlled case for the aforementioned standardized drive cycles is summarized in Table \ref{tab:standard_cycle}. We note that, the optimal PT controller shows improvement in terms of energy efficiency compared to their baseline counterpart in all the standardized drive cycles considered here. The most fuel consumption benefit is obtained for the UDDS drive cycle, and the least for HWFET.
\begin{table}[!htb]
\fontsize{8}{10}\selectfont
\centering
\caption{PT controller validation for standardized drive cycles.}\label{tab:standard_cycle}
\begin{tabular}{| L{0.2\columnwidth-2\tabcolsep-1.2\arrayrulewidth} | L{0.2\columnwidth-2\tabcolsep-1.2\arrayrulewidth} | L{0.2\columnwidth-2\tabcolsep-1.2\arrayrulewidth} | L{0.2\columnwidth-2\tabcolsep-1.2\arrayrulewidth} | 
L{0.2\columnwidth-2\tabcolsep-1.2\arrayrulewidth} | }
\hline
           \textbf{Drive Cycle [miles]} & \textbf{US06 \quad (8.0 )}         & \textbf{UDDS \quad(7.5)}         & \textbf{HWFET (10.3)}  & \textbf{Combined (25.7)}                       \\ \hline
\textbf{Baseline [MPGe]} & 26.4      & 28.2  & 32.5       &    29.1    \\  \hline
\textbf{PT Controller [MPGe]} & 27.8       & 30.4   & 38.1 & 31.8             \\  \hline
\textbf{Improvement [\%]} & 7.7     &  17.1    & 5.3 & 9.1             \\  \hline
	\end{tabular}
	\par
   \vspace{-0.15\skip\footins}
   \renewcommand{\footnoterule}{}
\end{table}

\begin{figure}[!ht]
	\centering
	\includegraphics[width=3.5 in]{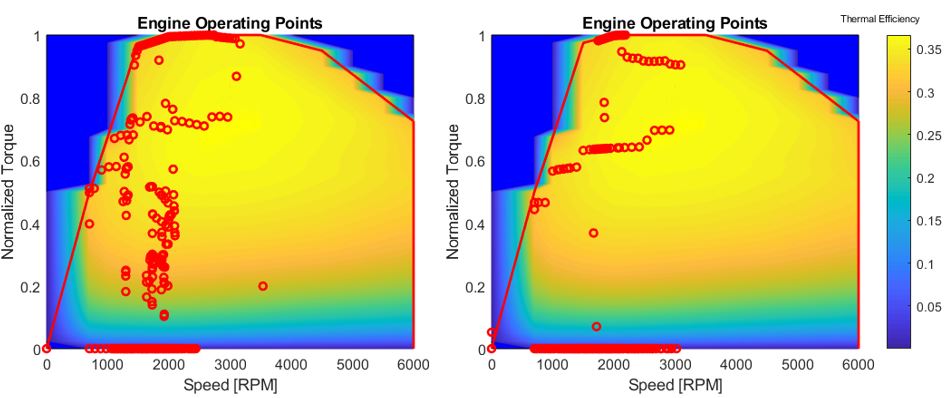}\caption{Engine operating points without and with the PT controller.}%
	\label{fig:operating-points}%
\end{figure}

In Fig. \ref{fig:operating-points}, we observe that the PT controller operates the engine at the most efficient brake specific fuel consumption regimes. In the baseline scenario, in contrast, there is a spread of operating points in non-efficient regimes. As the Pareto control policy yields online equivalent solution to DP, the benefits of the PT controller would be apparent for any heuristic approach used in the baseline scenario.

\subsection{VD Controller Performance}
\begin{figure}[ht]
	\centering
	\includegraphics[width=3.5in]{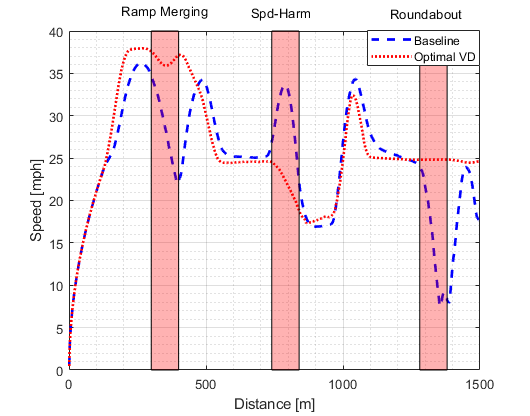}\caption{Average vehicle speed trajectories for high traffic volume.}%
	\label{fig:high-traffic}%
\end{figure}

Fig. \ref{fig:high-traffic} illustrates the average of vehicle speed profiles corresponding to PHEVs with the baseline scenario, and average of speed profiles corresponding to the CA-PHEVs with the VD controller travelling through the corridor under high traffic scenario. 
Note that the average speed profile of the baseline scenario in Fig. \ref{fig:high-traffic} shows that CAVs cruise with low speed in the control zones at ramp-merging and roundabout. The control zones of on-ramp merging and roundabout denote the upstream area of the entry to these bottlenecks. Human-driven vehicles on the ramp have to yield to the incoming vehicles from the main road. The human-driven vehicles exhibit stop-and-go behavior if the main road is very congested. Therefore, the resulting average speed in the baseline scenario (Fig. \ref{fig:high-traffic}) is low inside the control zones. On the other hand, with the VD controller, the CAVs can space themselves in such a way that they can enter the conflict zones without any stop-and-go driving behavior. As a result, the average speed of all CAVs is higher than the baseline in the control zones. Another interesting observation is that the average speed profile in the baseline scenario (Fig. \ref{fig:high-traffic}) is high at the SRZ control zone marked as \textit{Spd-Harm}. The SRZ is located at the end of a straight segment of Mcity (see Fig. \ref{fig:corridor}) which allows the vehicles to pick up speed. Since the conventional human driven vehicles do not have any information regarding the upcoming SRZ or the slower moving vehicle inside the SRZ, they pick up speed at the upstream of the SRZ and suddenly start decelerating at the entry of the SRZ. This behavior results in the backward propagating traffic wave. On the other hand, the CAVs know beforehand the state of the previous vehicles approaching the SRZ. Therefore, the CAVs adjust their speed inside the control zone of the SRZ in such a way that they have smooth entry at the SRZ while completely eliminating the backward propagating traffic wave.
We observe that under the optimal scenario of the VD controller, the average vehicle speed is more streamlined at the conflict zones compare to the average speed of the vehicles under baseline scenario. where speed oscillation upstream and downstream of the conflict zone is observed. The streamlined speed profile of the VD controller indicates a reduction in transient engine operation, resulting in 21.3\% better fuel efficiency compared to the baseline for high traffic volume. 

\begin{figure}[!ht]
	\centering
	\includegraphics[width=3.5 in]{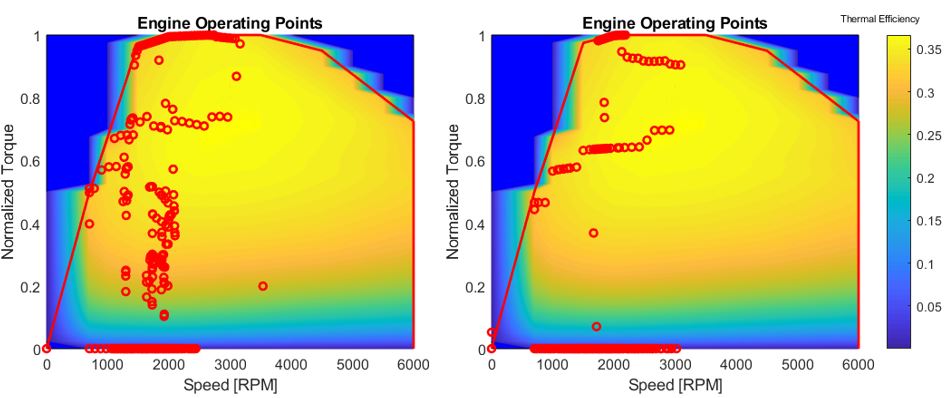}\caption{Engine operating points without and with the PT controller.}%
	\label{fig:operating-points}%
\end{figure}

\subsection{PT and VD Controller Performance}
To further investigate and quantify the individual as well as the combined contribution of the VD and PT controller, we consider four different simulation cases where different combination of the VD and PT controllers are active.
\begin{enumerate}
    \item \textbf{VD and PT controller inactive (Baseline):} None of the VD and PT controller are active in the CA-PHEV. The vehicle traverse the whole corridor by following the conventional traffic laws, and uses the Wiedemann \cite{wiedemann1974} car-following model. The engine and motor operating points are also determined by Audi A3's factory powertrain setting as determined by the baseline VESIM model.
    \item \textbf{VD controller active:} The CA-PHEVs traverse the route with optimal VD controller by communicating with other CA-PHEVs, but uses the Audi A3's factory powertrain model for selecting its torque set-point.
    \item \textbf{PT controller active:} The CA-PHEVs use the optimal PT controller to select the torque distribution between its engine and the IMG unit to meet the driver's demand. However, the vehicle is driven by conventional human drivers without the help of connectivity and automation.
    \item \textbf{VD and PT controller active:} The CA-PHEVs employ VD controller within the CZ to optimize its speed profile. For each optimal desired speed, the vehicle demands torque from the powertrain. The vehicle then uses the optimal PT controller to select the optimal power split between the engine and the IMG unit to meet the torque demanded.
\end{enumerate}

To evaluate the robustness of the these control cases in different traffic congestion, we consider low, medium and high traffic condition. \ref{tab:controller_cases} summarizes the VD and PT controller's impact on fuel efficiency in three different traffic condition. We observe significant increase in fuel efficiency in all traffic scenarios considered. Note that, the VD controller shows comparatively better improvement for high traffic volume as it streamlines the extreme stop-and-go driving behavior associated with high congestion scenario. On the other hand, the PT controller performs comparatively better in low traffic scenario. Although the VD and PT controller show energy improvement individually, their combination manages to obtain the most benefit in all traffic scenarios.

\begin{table}[!htb]
\fontsize{8}{10}\selectfont
\centering
\caption{Fuel Efficiency Improvement for VD and PT controller.}\label{tab:controller_cases}
\begin{tabular}{| L{0.4\columnwidth-2\tabcolsep-1.2\arrayrulewidth} | L{0.2\columnwidth-2\tabcolsep-1.2\arrayrulewidth} | L{0.2\columnwidth-2\tabcolsep-1.2\arrayrulewidth} | L{0.2\columnwidth-2\tabcolsep-1.2\arrayrulewidth} | 
}
\hline
           \textbf{Controllers / Traffic Level} & \textbf{Low}         & \textbf{Medium}         & \textbf{High}                      \\ \hline
\textbf{VD Controller Only [\%]} & 17.3      & 17.7   & 21.3           \\  \hline
\textbf{PT Controller Only [\%]} & 25.9     &  25.4    & 21.8          \\  \hline
\textbf{PT \& VD Controller Combined [\%]} & 34.0    &  32.7   & 29.2             \\  \hline
	\end{tabular}
	\par
   \vspace{-0.15\skip\footins}
   \renewcommand{\footnoterule}{}
\end{table}

\begin{figure}[ht]
	\centering
	\includegraphics[width=3.5in]{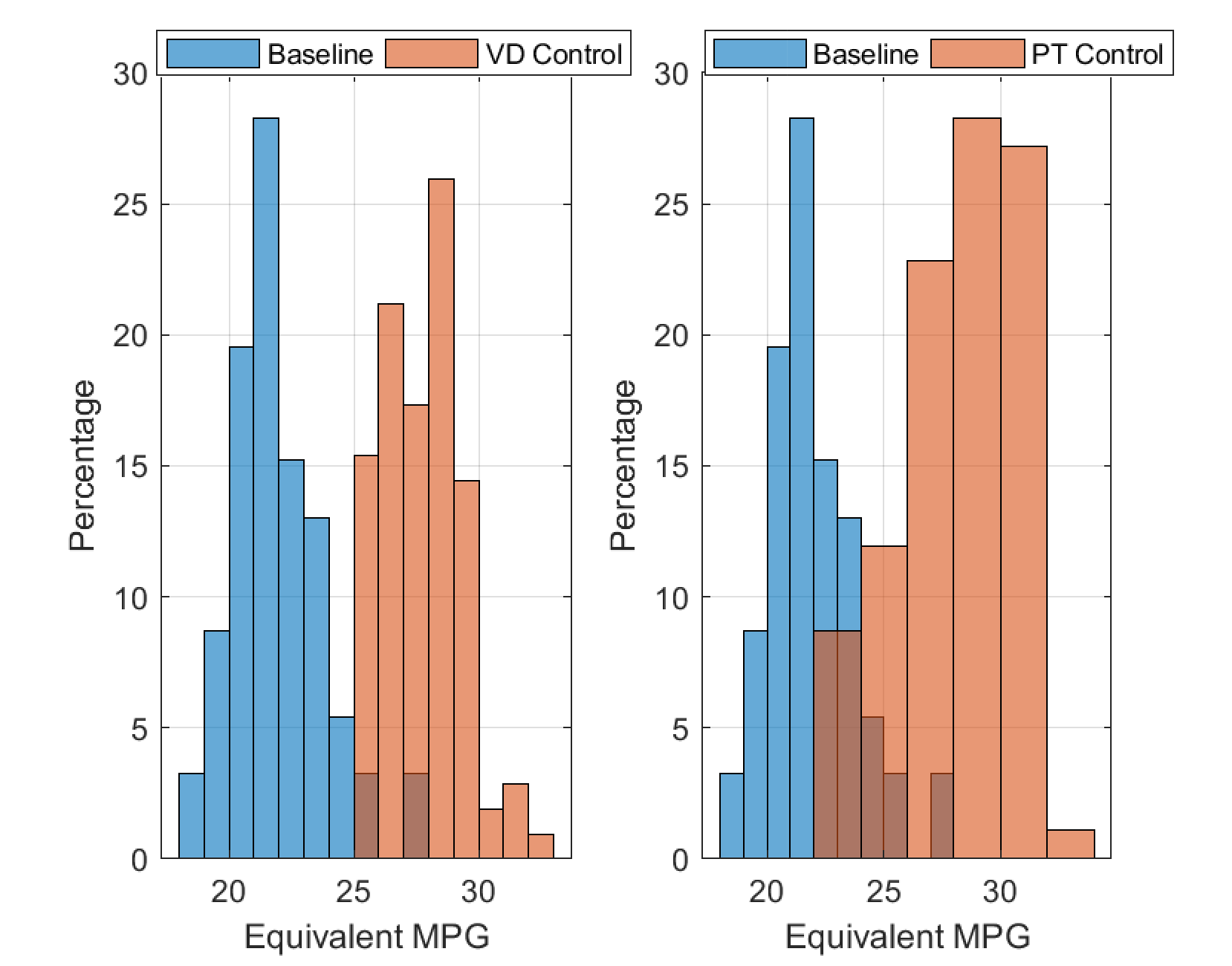}\caption{MPGe distribution in high traffic scenario.}%
	\label{fig:hist1}%
\end{figure}
\begin{figure}[ht]
	\centering
	\includegraphics[width=3.5in]{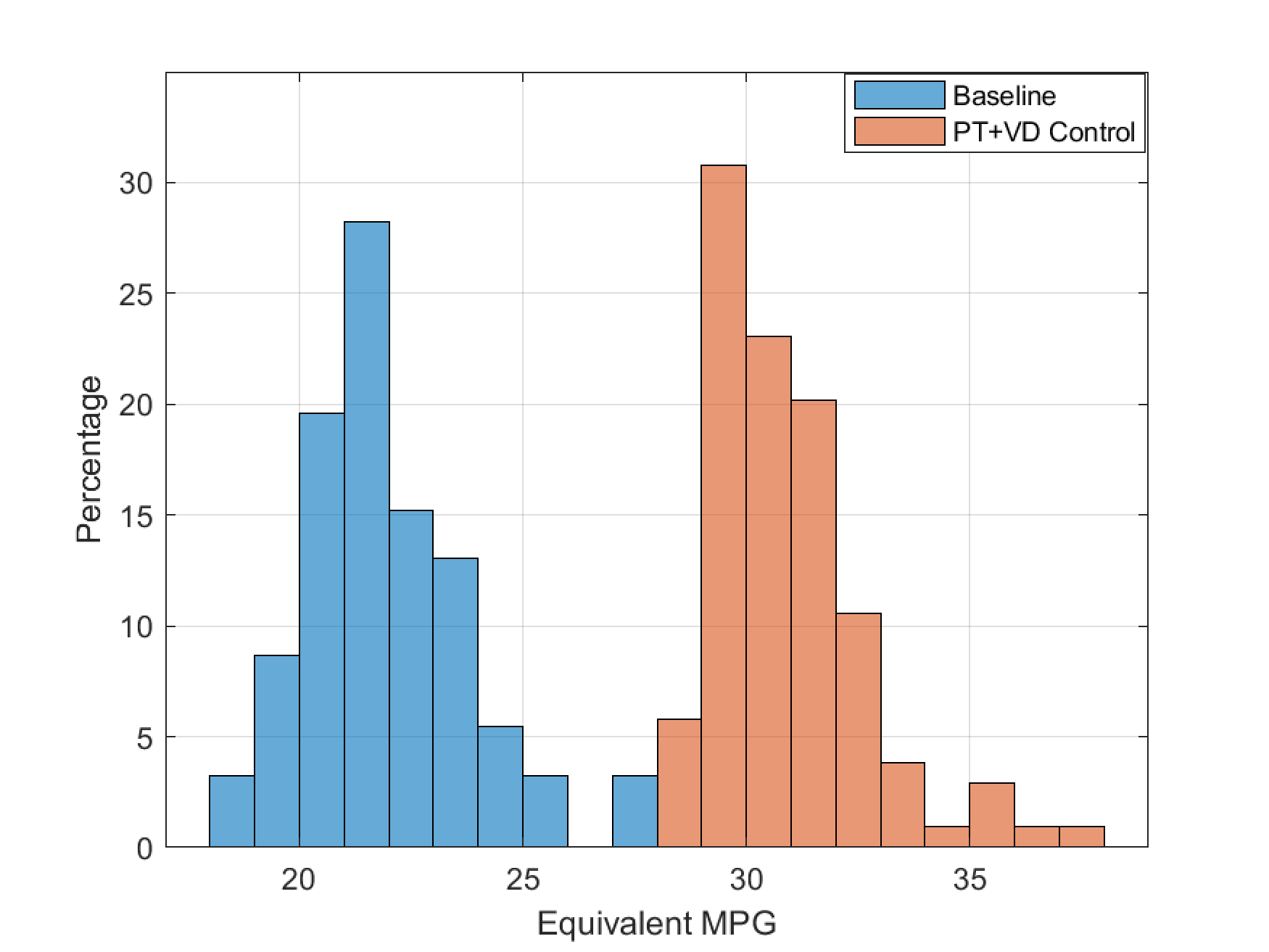}\caption{MPGe distribution in high traffic scenario.}%
	\label{fig:hist2}%
\end{figure}

\begin{figure}[ht]
	\centering
	\includegraphics[width=3.5in]{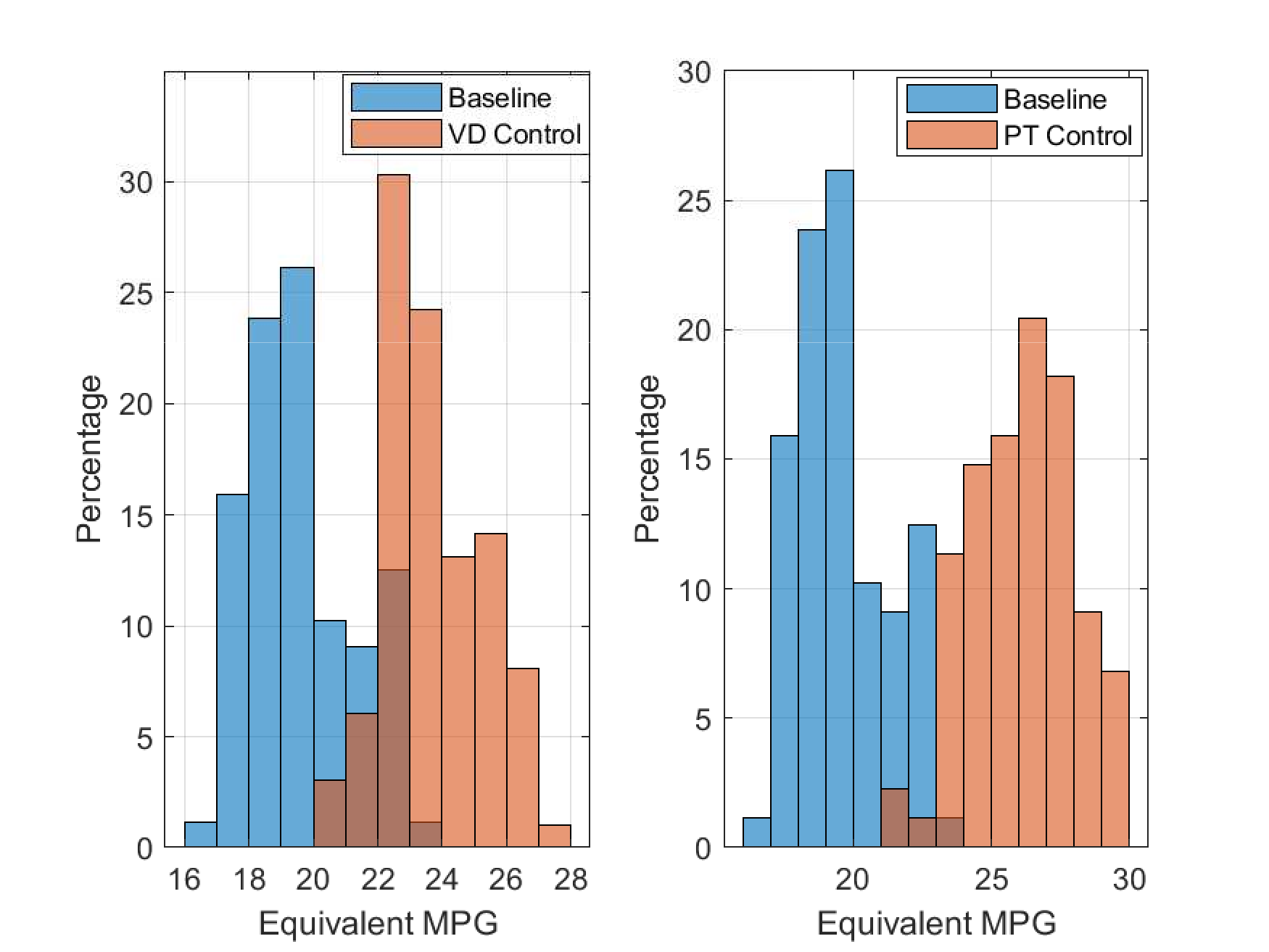}\caption{MPGe distribution in medium traffic scenario.}%
	\label{fig:hist3}%
\end{figure}
\begin{figure}[ht]
	\centering
	\includegraphics[width=3.5in]{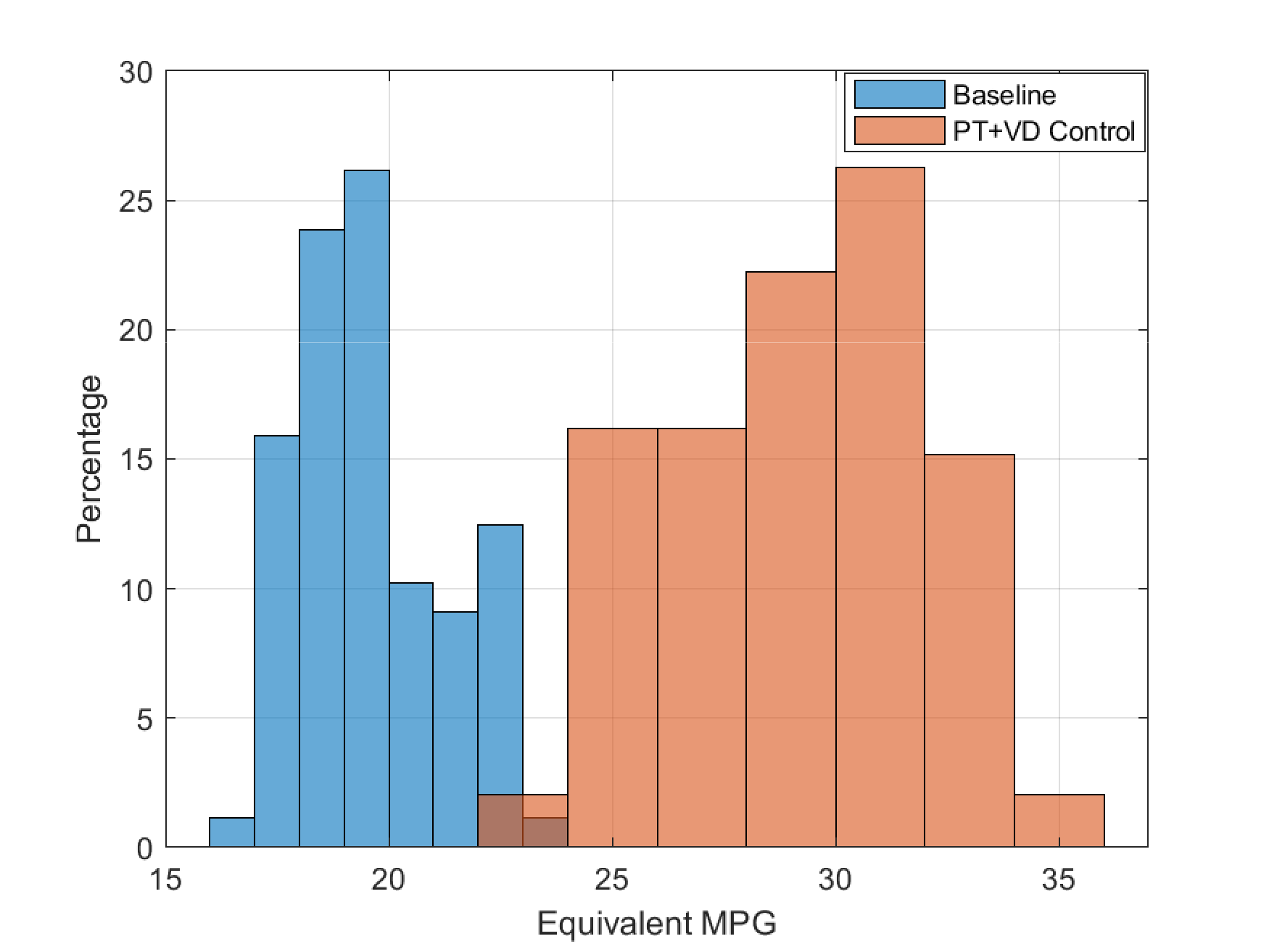}\caption{MPGe distribution in medium traffic scenario.}%
	\label{fig:hist4}%
\end{figure}

\begin{figure}[]
	\centering
	\includegraphics[width=3.5in]{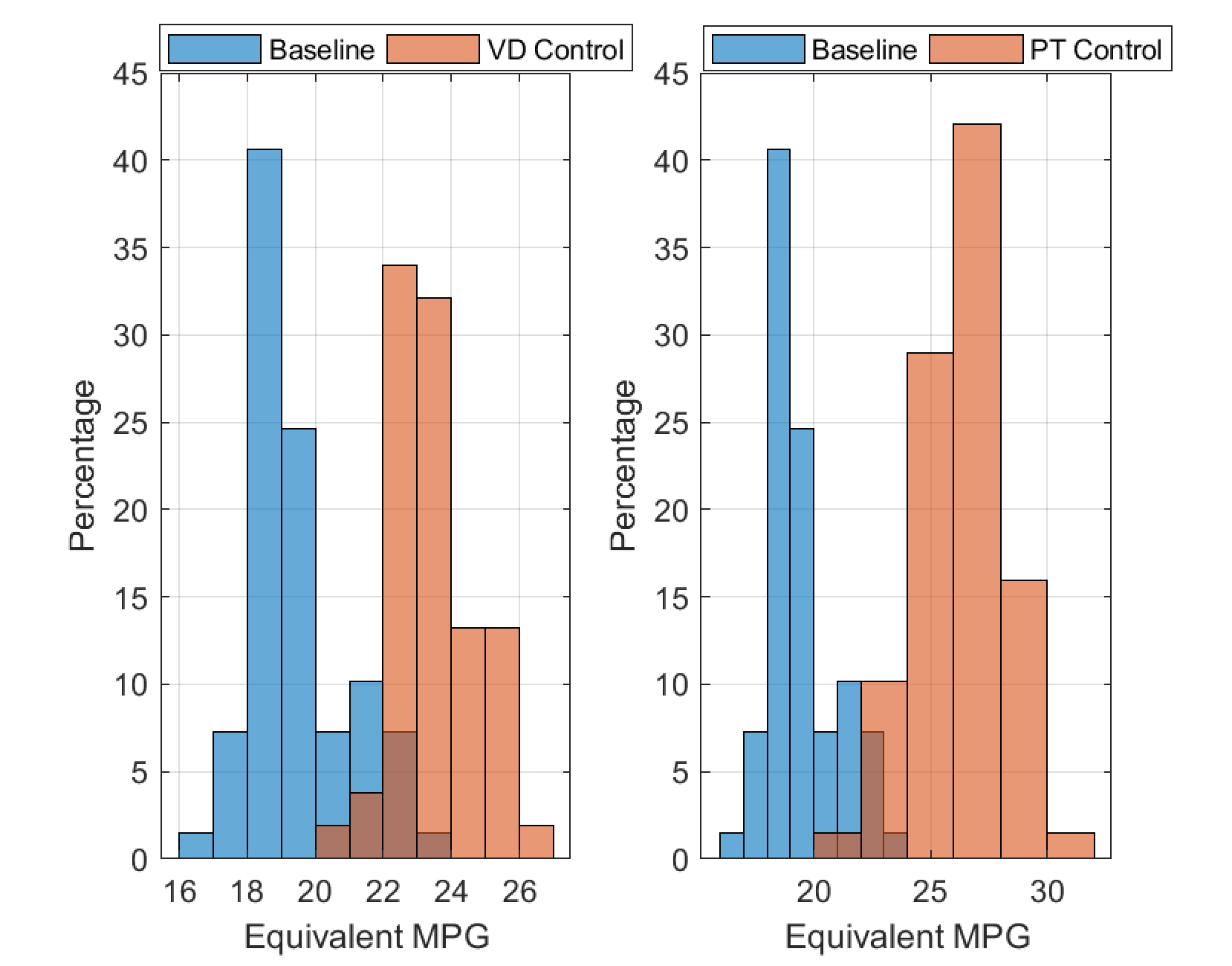}\caption{MPGe distribution in low traffic scenario.}%
	\label{fig:hist5}%
\end{figure}
\begin{figure}[]
	\centering
	\includegraphics[width=3.5in]{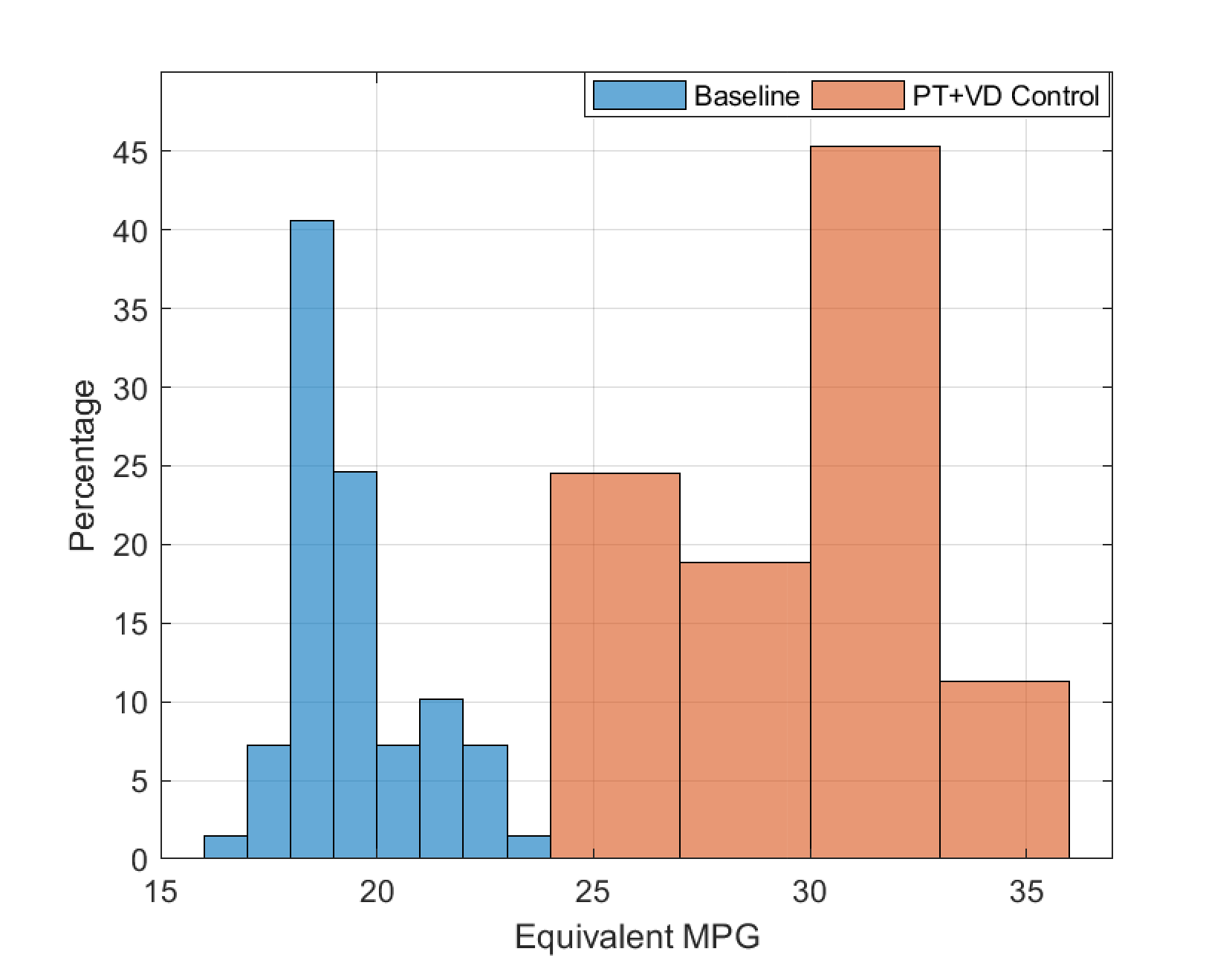}\caption{MPGe distribution in low traffic scenario.}%
	\label{fig:hist6}%
\end{figure}

The frequency distribution in terms of the vehicle MPGe in the four different control cases are illustrated in Figs. (\ref{fig:hist1})-(\ref{fig:hist2}) for high traffic volume, in Figs. (\ref{fig:hist3})-(\ref{fig:hist4}) for medium traffic volume, and in Figs. (\ref{fig:hist5})-(\ref{fig:hist6}) for low traffic volume. We analyze the frequency distributions presented in Figs. (\ref{fig:hist1})-(\ref{fig:hist6}) based on the skew property of distribution to get an insight about the nature and relative contribution of each controller cases. We observe that, for all the traffic volume cases, the baseline distribution has positive skew. The distribution resulting from the application of the VD controller is shifted to the right, but has the positive skew as well. However, when the PT controller is applied to the baseline cases, the positive skewed distribution is transformed into a negative skewed one. As a consequence, the implementation of the combined PT and VD controller gives rise to higher spread of the frequency distribution without showing any predominant skew. Note that, in the case of high traffic volume, the combined PT and VD distribution has positive skew. This is due to the fact that, the VD controller yields the most energy benefit in high traffic volume compare to the other traffic volume cases, thus transferring its skewed nature to the combined PT and VD controller effect.

\begin{table}[]
\fontsize{8}{10}\selectfont
\centering
\caption{Mean and standard deviation of the MPGe distribution for different traffic flow scenarios.}\label{tab:avg-std distribution}

\begin{tabular}{| L{0.3\columnwidth-2\tabcolsep-1.2\arrayrulewidth} | L{0.15\columnwidth-2\tabcolsep-1.2\arrayrulewidth} | L{0.15\columnwidth-2\tabcolsep-1.2\arrayrulewidth} | L{0.2\columnwidth-2\tabcolsep-1.2\arrayrulewidth} | 
L{0.2\columnwidth-2\tabcolsep-1.2\arrayrulewidth} | 
}
\hline
\textbf{}     & \textbf{Baseline} & \textbf{VD Active} & \textbf{PT \quad Active} & \textbf{VD+PT Active} \\ \hline
\textbf{}     & \multicolumn{4}{c|}{\textbf{Low Traffic}}                                \\ \hline
\textbf{Mean [MPGe]} & 19.5              & 23.5            & 26.3            & 29.8          \\ \hline
\textbf{Std. Deviation} & 1.5               & 1.3             & 1.9             & 2.8            \\ \hline
\textbf{}     & \multicolumn{4}{c|}{\textbf{Medium Traffic}}                             \\ \hline
\textbf{Mean [MPGe]} & 19.6              & 23.7             & 26.2            & 29.2           \\ \hline
\textbf{Std. Deviation} & 1.7              & 1.5              & 1.8             & 2.9           \\ \hline
\textbf{}     & \multicolumn{4}{c|}{\textbf{High Traffic}}                               \\ \hline
\textbf{Mean [MPGe]} & 21.9              & 27.8            & 28.1             & 30.9           \\ \hline
\textbf{Std. Deviation} & 1.8              & 1.51             & 2.5              & 1.7           \\ \hline
\end{tabular}
	\par
   \vspace{-0.15\skip\footins}
   \renewcommand{\footnoterule}{}
\end{table}
We further quantify the characteristics of the MPGe distribution under different control cases illustrated by Figs. (\ref{fig:hist1})-(\ref{fig:hist6}). In Table \ref{tab:avg-std distribution}, we summarize the mean and standard deviation of the distribution of the control cases under different traffic volumes. We observe that, for all the traffic volume cases, the VD controller reduces the standard deviation and increases the mean of the MPGe distribution compared to the baseline scenario. This implies that, by streamlining the vehicle speed profiles, the VD controller enables a more closely packed distribution of fuel consumption. On the other hand, the PT controller only case shows energy improvement compared to the baseline case for all traffic volumes as evident by the higher mean of the distribution, but shows increase in the standard deviation. Finally, the combined effect of the VD and PT controller is observed for the increase in the mean MPGe. However, the standard deviation of the distribution increases for both the low and medium traffic volume cases. Note that, for low and medium traffic volume, the increase in average MPGe for the combined VD and PT controlled case compared to the baseline case is additive in nature, i.e., the MPGe increase for only the case of the VD controller, and the MPGe increase for only the PT controlled case can be added to obtain the MPGe increase for their combined control. The high traffic volume case, however, slightly varies from this observation.

\section{Concluding Remarks}
In this paper, we proposed a two-level control architecture with the aim to optimize simultaneously the vehicle-level and powertrain-level operation of a PHEV. With the proposed approach, we can optimize both the speed profile of a vehicle by eliminating stop-and-go driving and the powertrain efficiency.
We applied the proposed architecture to the operation of CA-PHEVs over a range of real-world driving scenarios deemed characteristic of typical commute that included a merging roadway, a speed reduction zone, and a roundabout. Ongoing work considers additional scenarios, including intersections, while incorporating the state and control constraints in the analytical solution of the VD controller. Mixed traffic scenario with interacting CAVs and human driven vehicles is also under consideration.

While the potential benefits of full penetration of CAVs to alleviate traffic congestion and reduce fuel consumption have become apparent, different penetrations of CAVs can alter significantly the efficiency of the entire system. 
For example, one particular question that still remains unanswered is ``what is the minimum number of vehicles that need to be connected so that to start realizing the potential benefits?'' The assumption of perfect information seems to impose barriers in a potential implementation and deployment of the proposed framework. Although it is relatively straightforward to extend our results in the case that this assumption is relaxed, future research should investigate the implications of having information with errors and/or delays to the system behavior. Finally, considering lane changing with a diverse set of CAVs would eventually aim at addressing the remaining practical consequences of implementing this framework.

\bibliographystyle{ieeetr} 
\bibliography{acc_pt_vd_ref}

\begin{thebibliography}{10}

\bibitem{Malikopoulos2014b}
A.~A. Malikopoulos, ``{Supervisory Power Management Control Algorithms for
  Hybrid Electric Vehicles: A Survey},'' {\em IEEE Transactions on Intelligent
  Transportation Systems}, vol.~15, no.~5, pp.~1869--1885, 2014.

\bibitem{Malikopoulos2008}
A.~A. Malikopoulos, {\em {Real-Time, Self-Learning Identification and
  Stochastic Optimal Control of Advanced Powertrain Systems}}.
\newblock ProQuest, September 2011.

\bibitem{Malikopoulos2013}
A.~A. Malikopoulos and J.~P. Aguilar, ``{An Optimization Framework for Driver
  Feedback Systems},'' {\em IEEE Transactions on Intelligent Transportation
  Systems}, vol.~14, no.~2, pp.~955--964, 2013.

\bibitem{Athans1969}
M.~Athans, ``{A unified approach to the vehicle-merging problem},'' {\em
  Transportation Research}, vol.~3, no.~1, pp.~123--133, 1969.

\bibitem{Varaiya1993}
P.~Varaiya, ``Smart cars on smart roads: problems of control,'' {\em IEEE
  Transactions on Automatic Control}, vol.~38, no.~2, pp.~195--207, 1993.

\bibitem{Dresner2004}
K.~Dresner and P.~Stone, ``{Multiagent traffic management: a reservation-based
  intersection control mechanism},'' in {\em Proceedings of the Third
  International Joint Conference on Autonomous Agents and Multiagents Systems},
  pp.~530--537, 2004.

\bibitem{Dresner2008}
K.~Dresner and P.~Stone, ``A multiagent approach to autonomous intersection
  management,'' {\em Journal of artificial intelligence research}, vol.~31,
  pp.~591--656, 2008.

\bibitem{DeLaFortelle2010}
A.~{de La Fortelle}, ``{Analysis of reservation algorithms for cooperative
  planning at intersections},'' {\em 13th International IEEE Conference on
  Intelligent Transportation Systems}, pp.~445--449, Sept. 2010.

\bibitem{Huang2012}
S.~Huang, A.~Sadek, and Y.~Zhao, ``{Assessing the Mobility and Environmental
  Benefits of Reservation-Based Intelligent Intersections Using an Integrated
  Simulator},'' {\em IEEE Transactions on Intelligent Transportation Systems},
  vol.~13, no.~3, pp.~1201--1214, 2012.

\bibitem{Zohdy2012}
I.~H. Zohdy, R.~K. Kamalanathsharma, and H.~Rakha, ``{Intersection management
  for autonomous vehicles using iCACC},'' {\em 2012 15th International IEEE
  Conference on Intelligent Transportation Systems}, pp.~1109--1114, 2012.

\bibitem{Yan2009}
F.~Yan, M.~Dridi, and A.~{El Moudni}, ``{Autonomous vehicle sequencing
  algorithm at isolated intersections},'' {\em 2009 12th International IEEE
  Conference on Intelligent Transportation Systems}, pp.~1--6, 2009.

\bibitem{kim2014}
K.-D. Kim and P.~Kumar, ``{An MPC-Based Approach to Provable System-Wide Safety
  and Liveness of Autonomous Ground Traffic},'' {\em IEEE Transactions on
  Automatic Control}, vol.~59, no.~12, pp.~3341--3356, 2014.

\bibitem{Rios-Torres2}
J.~Rios-Torres and A.~A. Malikopoulos, ``{Automated and Cooperative Vehicle
  Merging at Highway On-Ramps},'' {\em IEEE Transactions on Intelligent
  Transportation Systems}, vol.~18, no.~4, pp.~780--789, 2017.

\bibitem{Ntousakis:2016aa}
I.~A. Ntousakis, I.~K. Nikolos, and M.~Papageorgiou, ``Optimal vehicle
  trajectory planning in the context of cooperative merging on highways,'' {\em
  Transportation Research Part C: Emerging Technologies}, vol.~71,
  pp.~464--488, 2016.

\bibitem{Malikopoulos2017}
A.~A. Malikopoulos, C.~G. Cassandras, and Y.~J. Zhang, ``A decentralized
  energy-optimal control framework for connected automated vehicles at
  signal-free intersections,'' {\em Automatica}, vol.~93, pp.~244 -- 256, 2018.

\bibitem{Mahbub2019ACC}
A.~M.~I. Mahbub, L.~Zhao, D.~Assanis, and A.~A. Malikopoulos, ``{Energy-Optimal
  Coordination of Connected and Automated Vehicles at Multiple
  Intersections},'' in {\em Proceedings of 2019 American Control Conference},
  pp.~2664--2669, 2019.

\bibitem{Malikopoulos2018a}
L.~Zhao, A.~A. Malikopoulos, and J.~Rios-Torres, ``Optimal control of connected
  and automated vehicles at roundabouts: An investigation in a mixed-traffic
  environment,'' in {\em 15th IFAC Symposium on Control in Transportation
  Systems}, pp.~73--78, 2018.

\bibitem{Malikopoulos2018c}
A.~A. Malikopoulos, S.~Hong, B.~Park, J.~Lee, and S.~Ryu, ``Optimal control for
  speed harmonization of automated vehicles,'' {\em IEEE Transactions on
  Intelligent Transportation Systems}, vol.~20, no.~7, pp.~2405--2417, 2018.

\bibitem{Zhao2018}
L.~Zhao and A.~A. Malikopoulos, ``Decentralized optimal control of connected
  and automated vehicles in a corridor,'' in {\em 2018 21st International
  Conference on Intelligent Transportation Systems (ITSC)}, pp.~1252--1257, Nov
  2018.

\bibitem{Mahbub2019CDC}
A.~M.~I. Mahbub and A.~A. Malikopoulos, ``Conditions for state and control
  constraint activation in coordination of connected and automated vehicles,''
  {\em Proceedings of 2020 American Control Conference}, 2020 (to appear) arXiv
  preprint arxiv:1903.11189.

\bibitem{beaver2019demonstration}
L.~E. Beaver, B.~Chalaki, A.~M.~I. Mahbub, L.~Zhao, R.~Zayas, and A.~A.
  Malikopoulos, ``Demonstration of a time-efficient mobility system using a
  scaled smart city,'' {\em Vehicle System Dynamics}, 2019 (to appear) arXiv
  preprint arxiv:1903.01632.

\bibitem{Malikopoulos2016a}
J.~Rios-Torres and A.~A. Malikopoulos, ``{A Survey on Coordination of Connected
  and Automated Vehicles at Intersections and Merging at Highway On-Ramps},''
  {\em IEEE Transactions on Intelligent Transportation Systems}, vol.~18,
  no.~5, pp.~1066--1077, 2017.

\bibitem{Malikopoulos2016AMO}
A.~A. Malikopoulos, ``A multiobjective optimization framework for online
  stochastic optimal control in hybrid electric vehicles,'' {\em IEEE
  Transactions on Control Systems Technology}, vol.~24, pp.~440--450, 2016.

\bibitem{Malikopoulos2014c}
M.~Shaltout, A.~A. Malikopoulos, S.~Pannala, and D.~Chen, ``A consumer-oriented
  control framework for performance analysis in hybrid electric vehicles,''
  {\em IEEE Transactions on Control Systems Technology}, vol.~23, no.~4,
  pp.~1451--1464, 2015.

\bibitem{Malikopoulos2013a}
A.~A. Malikopoulos, ``Stochastic optimal control for series hybrid electric
  vehicles,'' in {\em Proceedings of the 2013 American Control Conference},
  pp.~1189--1194, 2013.

\bibitem{Malikopoulos2015_ITS_HEV}
A.~A. Malikopoulos, ``Pareto efficient policy for supervisory power management
  control,'' in {\em Proceedings of the 2015 IEEE 18th International Conference
  on Intelligent Transportation Systems}, September 15-18, 2015.

\bibitem{Lin2003}
C.-C. Lin, H.~Peng, J.~W. Grizzle, and J.-m. Kang, ``{Power Management Strategy
  for a Parallel Hybrid Electric Truck},'' {\em IEEE Transactions on Control
  Systems Technology}, vol.~11, no.~6, pp.~839--849, 2003.

\bibitem{Lin2004}
C.-C. Lin, H.~Peng, and J.~W. Grizzle, ``{A stochastic control strategy for
  hybrid electric vehicles},'' in {\em Proceedings of the 2004 American Control
  Conference}, vol.~5, pp.~4710--4715, 2004.

\bibitem{Tate2010}
E.~D. Tate, J.~W. Grizzle, and H.~Peng, ``{SP-SDP for Fuel Consumption and
  Tailpipe Emissions Minimization in an EVT Hybrid},'' {\em IEEE Transactions
  on Control Systems Technology}, vol.~18, no.~3, pp.~1--16, 2010.

\bibitem{Tate2007}
E.~D.~J. Tate, J.~W. Grizzle, and H.~Peng, ``{Shortest path stochastic control
  for hybrid electric vehicles},'' {\em International Journal of Robust and
  Nonlinear Control}, vol.~18, pp.~1409--1429, 2008.

\bibitem{Opila2008}
D.~F. Opila, D.~Aswani, R.~McGee, J.~a. Cook, and J.~W. Grizzle,
  ``{Incorporating drivability metrics into optimal energy management
  strategies for Hybrid Vehicles},'' {\em 2008 47th IEEE Conference on Decision
  and Control}, pp.~4382--4389, 2008.

\bibitem{Paganelli2001}
G.~Paganelli, M.~Tateno, A.~Brahma, G.~Rizzoni, and Y.~Guezennec, ``{Control
  development for a hybrid-electric sport-utility vehicle: strategy,
  implementation and field test results},'' in {\em Proceedings of the 2001
  American Control Conference}, vol.~6, pp.~5064--5069, 2001.

\bibitem{Sciarretta2004}
A.~Sciarretta, M.~Back, and L.~Guzzella, ``{Optimal Control of Parallel Hybrid
  Electric Vehicles},'' {\em IEEE Transactions on Control Systems Technology},
  vol.~12, no.~3, pp.~352--363, 2004.

\bibitem{Musardo2005}
C.~Musardo, G.~Rizzoni, and B.~Staccia, ``{A-ECMS: An Adaptive Algorithm for
  Hybrid Electric Vehicle Energy Management},'' in {\em Proceedings of the 44th
  IEEE Conference on Decision and Control, and the European Control
  Conference}, pp.~1816--1823, 2005.

\bibitem{Pisu2007}
P.~Pisu and G.~Rizzoni, ``{A Comparative Study of Supervisory Control
  Strategies for Hybrid Electric Vehicles},'' {\em IEEE Transactions on Control
  Systems Technology}, vol.~15, no.~3, pp.~506--518, 2007.

\bibitem{Borhan2012}
H.~Borhan, A.~Vahidi, A.~M. Phillips, M.~L. Kuang, I.~V. Kolmanovsky, and
  S.~{Di Cairano}, ``{MPC-Based Energy Management of a Power-Split Hybrid
  Electric Vehicle},'' {\em IEEE Transactions on Control Systems Technology},
  vol.~20, no.~3, pp.~593--603, 2012.

\bibitem{Johannesson2007}
L.~Johannesson, M.~\AA~sbog\aa rd, and B.~Egardt, ``{Assessing the Potential of
  Predictive Control for Hybrid Vehicle Powertrains Using Stochastic Dynamic
  Programming},'' vol.~8, no.~1, pp.~71--83, 2007.

\bibitem{Ambuhl2009}
D.~Ambuhl and L.~Guzzella, ``{Predictive Reference Signal Generator for Hybrid
  Electric Vehicles},'' {\em IEEE Transactions on Vehicular Technology},
  vol.~58, no.~9, pp.~4730--4740, 2009.

\bibitem{Malikopoulos2012}
A.~A. Malikopoulos and J.~P. Aguilar, ``{Optimization of driving styles for
  fuel economy improvement},'' in {\em Proceedings of the 2012 15th
  International IEEE Conference on Intelligent Transportation Systems},
  pp.~194,199, 2012.

\bibitem{Malikopoulos2013d}
A.~A. Malikopoulos and J.~P. Aguilar, ``An optimization framework for driver
  feedback systems,'' {\em IEEE Transactions on Intelligent Transportation
  Systems}, vol.~14, no.~2, pp.~955--964, 2013.

\bibitem{Huang2011}
X.~Huang, Y.~Tan, and X.~He, ``An intelligent multifeature statistical approach
  for the discrimination of driving conditions of a hybrid electric vehicle,''
  {\em IEEE Transactions on Intelligent Transportation Systems}, vol.~12,
  no.~2, pp.~453--465, 2011.

\bibitem{Luo2015a}
Y.~Luo, T.~Chen, and K.~Li, ``Multi-objective decoupling algorithm for active
  distance control of intelligent hybrid electric vehicle,'' {\em Mechanical
  Systems and Signal Processing}, vol.~64-65, pp.~29--45, 2015.

\bibitem{Luo2015b}
Y.~Luo, T.~Chen, S.~Zhang, and K.~Li, ``Intelligent hybrid electric vehicle acc
  with coordinated control of tracking ability, fuel economy, and ride
  comfort,'' {\em IEEE Transactions on Intelligent Transportation Systems},
  vol.~16, pp.~2303--2308, Aug 2015.

\bibitem{Li2017a}
L.~Li, X.~Wang, and J.~Song, ``Fuel consumption optimization for smart hybrid
  electric vehicle during a car-following process,'' {\em Mechanical Systems
  and Signal Processing}, vol.~87, pp.~17--19, 2017.

\bibitem{Zhao2019CCTAa}
L.~Zhao, A.~Mahbub, and A.~A. Malikopoulos, ``Optimal vehicle dynamics and
  powertrain control for connected and automated vehicles,'' {\em Proceedings
  of 2019 IEEE Conference on Control Technology and Applications}, pp.~33--38,
  2019.

\bibitem{Malikopoulos2016b}
A.~A. Malikopoulos, ``A duality framework for stochastic optimal control of
  complex systems,'' {\em IEEE Transactions on Automatic Control}, vol.~18,
  no.~4, pp.~780--789, 2016.

\bibitem{ptv2014ptv}
P.~Group {\em et~al.}, ``Ptv vissim 7 user manual,'' {\em Germany: PTV GROUP},
  2014.

\bibitem{wiedemann1974}
R.~Wiedemann, ``Simulation des strassenverkehrsflusses.,'' 1974.

\bibitem{vesim-1000037126/PAPER}
C.-C. Lin, Z.~Filipi, Y.~Wang, L.~Louca, H.~Peng, D.~Assanis, and J.~Stein,
  ``{Integrated, Feed-Forward Hybrid Electric Vehicle Simulation in SIMULINK
  and its Use for Power Management Studies},'' 2001.

\end{thebibliography}

\section{Contact Information}
A M Ishtiaque Mahbub, \newline
mahbub@udel.edu \newline

Andreas A. Malikopoulos, \newline
andreas@udel.edu

\section{Acknowledgments}
This research was supported by ARPAE's NEXTCAR program under the award number DE-AR0000796.

\newpage
\section{Definitions, Acronyms, Abbreviations}

\begin{table}[h]
\centering
\begin{tabular}{L{0.1\textwidth} L{0.33\textwidth}}
\textbf{CA-PHEV} & Connected automated plug-in hybrid electric vehicle \\
\textbf{CAV} & Connected automated vehicle \\
\textbf{CZ} & Control zone \\
\textbf{DP} & Dynamic programming \\
\textbf{ECMS} & Equivalent consumption minimization strategy \\
\textbf{HWFET} & Highway fuel economy driving schedule \\
\textbf{HEV} & Hybrid electric vehicle \\
\textbf{MPGe} & Miles-per-gallon of gasoline equivalent \\
\textbf{PHEV} & Plugin hybrid electric vehicle \\
\textbf{PT} & Powertrain \\
\textbf{SRZ} & Speed reduction zone \\
\textbf{SOC} & State-of-charge \\
\textbf{UDDS} & Urban dynamometer driving schedule \\
\textbf{VD} & Vehicle dynamics \\
\textbf{V2V} & Vehicle to vehicle \\
\textbf{V2I} & Vehicle to infrastructure \\
\textbf{VPH} & Vehicle per hour \\
\textbf{VESIM} & Vehicle simulation \\
\end{tabular}
\end{table}

\end{document}